\documentclass[a4paper,reqno,10pt,oneside]{amsart}
\usepackage{amsmath}
\usepackage[left=2.61cm,right=2.61cm,top=2.72cm,bottom=2.72cm]{geometry}
\usepackage[colorlinks=true,urlcolor=blue,
citecolor=red,linkcolor=blue,linktocpage,pdfpagelabels,
bookmarksnumbered,bookmarksopen]{hyperref}
\usepackage[english]{babel}
\usepackage{graphicx}
\usepackage[toc,page,header]{appendix}
\usepackage{mathrsfs}
\usepackage{amssymb,amsmath, amsfonts, amsthm}
\usepackage{latexsym}
\usepackage{enumitem}

\allowdisplaybreaks[1]

\numberwithin{equation}{section}

\renewcommand{\(}{\left(}
\renewcommand{\)}{\right)}
\renewcommand{\[}{\left[}
\renewcommand{\]}{\right]}

\newtheorem{theorem}{Theorem}[section]
\newtheorem{proposition}[theorem]{Proposition}

\theoremstyle{definition}
\newtheorem{remark}[theorem]{Remark}
\theoremstyle{definition}
\newtheorem{definition}[theorem]{Definition}
\theoremstyle{definition}

\renewcommand{\le}{\leqslant}
\renewcommand{\ge}{\geqslant}

\newcommand{\N}{\mathbb{N}}
\renewcommand{\S}{\mathbb{S}}

\newcommand{\beq}{\begin{equation}}
\newcommand{\eeq}{\end{equation}}
\newcommand{\beqs}{\begin{equation*}}
\newcommand{\eeqs}{\end{equation*}}
\newcommand{\beqn}{\begin{eqnarray}}
\newcommand{\eeqn}{\end{eqnarray}}
\newcommand{\beqns}{\begin{eqnarray*}}
\newcommand{\eeqns}{\end{eqnarray*}}
\newcommand{\bdoc}{\begin{document}}
\newcommand{\edoc}{\end{document}}
\newcommand{\be}{\begin{enumerate}}
\newcommand{\ee}{\end{enumerate}}
\newcommand{\bdescr}{\begin{description}}
\newcommand{\edescr}{\end{description}}
\newcommand{\ba}{\begin{array}}
\newcommand{\ea}{\end{array}}
\newcommand{\intR}{\int_{\mathbb R^N}}
\newcommand{\R}{\mathbb R}
\newcommand{\RN}{\mathbb{R}^N}
\newcommand{\B}{\mathbb B}
\newcommand{\C}{\mathbb C}
\renewcommand{\H}{\mathcal H}
\renewcommand{\L}{\mathbb L}
\newcommand{\parallelsum}{\mathbin{\!/\mkern-5mu/\!}}
\newcommand{\e}{\varepsilon}
\newcommand{\SD}{\Sigma_D}
 \renewcommand{\(}{\left(}
\renewcommand{\)}{\right)}
\renewcommand{\[}{\left[}
\renewcommand{\]}{\right]}
\renewcommand{\appendixpagename}{\centering Appendix}

\newcommand{\todo}[1]{\text{\colorbox{yellow}{#1}}}

\begin{document}
\title[Symmetry breaking and instability for semilinear equations]{Symmetry breaking and instability for semilinear elliptic equations in spherical sectors and cones}

\author{Giulio Ciraolo}
\address{G. Ciraolo. Dipartimento di Matematica "Federigo Enriques",
Universit\`a degli Studi di Milano, Via Cesare Saldini 50, 20133 Milano, Italy}
\email{giulio.ciraolo@unimi.it}

\author{Filomena Pacella}
\address{F. Pacella. Dipartimento di Matematica, Sapienza Universit\`a di Roma,  P.le Aldo Moro 2, 00185 Roma, Italy}
\email{pacella@mat.uniroma1.it}

\author{Camilla Chiara Polvara}
\address{C.C. Polvara. Dipartimento di Matematica "Federigo Enriques",
Universit\`a degli Studi di Milano, Via Cesare Saldini 50, 20133 Milano, Italy}
\email{camilla.polvara@unimi.it}

\thanks{\emph{Acknowledgements.} Research partially supported by Gruppo Nazionale per l'Analisi Matematica, la Pro\-ba\-bi\-li\-t\`a e le loro Applicazioni (GNAMPA) of the Istituto Nazionale di Alta Matematica (INdAM)}

\begin{abstract}
We consider semilinear elliptic equations with mixed boundary conditions in spherical sectors inside a cone. The aim of the paper is to show that a radial symmetry result of Gidas-Ni-Nirenberg type for positive solutions does not hold in general nonconvex cones. This symmetry breaking result is achieved by studying the Morse index of radial positive solutions and analyzing how it depends on the domain D on the unit sphere which spans the cone.
In particular it is proved that the Neumann eigenvalues of the Laplace Beltrami operator on D play a role in computing the Morse index.
A similar breaking of symmetry result is obtained for the positive solutions of the critical Neumann problem in the whole unbounded cone. In this case it is proved that the standard bubbles, which are the only radial solutions, become unstable for a class of nonconvex cones. 
\end{abstract}

\maketitle

\section{Introduction}
Let $D$ be a smooth domain on the unit sphere $\S^{N-1}$, with $N\ge 3$, and let $\Sigma_D \subset \R^N$ be the cone spanned by $D$, namely
\beq\label{def:conespanned}
\Sigma_D:=\{x \in \R^{N}; \ x=s q,\ q\in D,\ s\in(0,+\infty)\}.
\eeq
Let us consider the spherical sector $S_D$ obtained by intersecting $\Sigma_D$ with the unit ball $B_1$
\beq\label{def:bddconespanned}
S_D:=\{x \in \R^{N}; \ x=s q,\ q\in D,\ s\in(0,1)\}.
\eeq
Then the relative boundary of $S_D$ (i.e.  the part of $\partial S_D$ which lies inside the cone ) is just $D$,  while the remaining part of $\partial S_D$ will be denoted by $\Gamma_1 $ and it is:
\beq
\Gamma_1 :=\partial \Sigma_D\cap \partial S_D \,.
\eeq
In the spherical sector $S_D$ we consider the following mixed boundary value problem:
\beq\label{eq:semilinearelliptic}
\begin{cases}
-\Delta u=f(u)& \text{in } S_D\\
\frac{\partial u}{\partial \nu}=0 & \text{on } \Gamma_1 \\
u=0& \text{on } D\\
\end{cases}
\eeq
where $f:\R \to \R $ is a locally Lipschitz continuous function and $\nu$ denotes the exterior unit normal vector.

We study positive weak solutions $u$ of \eqref{eq:semilinearelliptic} which are functions in the Sobolev space $ H_0^1(S_D\cup\Gamma_1 )$ which satisfy:
$$
\int_{S_D} (\nabla u\cdot\nabla\varphi-f(u)\varphi)dx=0
$$
for any test function $\varphi\in H_0^1(S_D\cup\Gamma_1 )$.  Here, $H_0^1(S_D\cup\Gamma_1 )$ is defined as the completion of $C^1_c(S_D\cup\Gamma_1 )$ with respect to the $H^1$ norm.

It is not difficult to see that, under suitable assumptions on the nonlinearity $f$, positive weak solutions of \eqref{eq:semilinearelliptic} exist. They can be easily obtained by adapting standard variational methods to the case of the spherical sector $S_D$.  In particular by considering the subspace $H_{0,rad}^1(S_D\cup \Gamma_1 )$ of the radial functions in $H_0^1(S_D\cup\Gamma_1 )$ we can obtain positive radial weak solutions of \eqref{eq:semilinearelliptic}.

The question we address in the present paper is whether or not all positive solutions of \eqref{eq:semilinearelliptic} are radial. 

It is well known that if $S_D$ is the unit ball $B_1$, i.e. if $D=\S^{N-1}$, then all positive classical solutions of \eqref{eq:semilinearelliptic} are radial by the famous theorem of  Gidas-Ni-Nirenberg \cite{GNN}.
This result was later extended to solutions in $W^{2,N}_{loc}(B_1)\cap C(\bar B_1)$ in the paper \cite{BN} and to weak solutions in $H_0^1(B_1)\cap L^\infty(B_1)$ in \cite{D}.

The method for proving symmetry of solutions used in these papers is the famous moving planes method which relies on different forms of maximum principles. It goes back to the classic paper of Serrin, \cite{Serrin}, and Alexandrov \cite{Alexandrov}, where it was introduced to get the radial symmetry of constant mean curvature surfaces and of domains admitting solutions of overdeterminated problems.

In the case of problem \eqref{eq:semilinearelliptic} the method of moving planes does not seem easily applicable due to the fact that a sector $S_D$ is not, in general, symmetric with respect to a suitable family of hyperplanes.  A quite involved modification of the moving planes method  was used in \cite{BP} to get the radial symmetry of classical $C^2$-solutions of \eqref{eq:semilinearelliptic} only in dimensions 2 and for spherical sectors with angle $\alpha\in (0,\pi)$, see also \cite{DePa} for the case of the unbounded cone. An extension of that proof to higher dimensions does not seem possible.

On the other side results analogous to Serrin and Alexandrov's ones have been recently obtained in the relative setting of cones, in all dimension $N\ge 3$ under the hypothesis that $\Sigma_D$ is a convex cone and for positive nonlinearities \cite{PT} (see also \cite{ChP} and \cite{CR2}). The proof of \cite{PT} is based on integral identities.

In analogy with the result of \cite{PT} it is natural to expect that radial symmetry for the positive solution of \eqref{eq:semilinearelliptic} should hold in any dimension $N\ge 3$ whenever the cone is convex. Recently this has been proved in \cite{DPPV} but only when the nonlinarity $f$ is nonnegative. The proof in \cite{DPPV} is based on integral identities as the one in \cite{PLL} and \cite{KP} and holds for equations involving more general operators. 

Giving the fact that all these symmetry results require the convexity of the cone it is interesting to understand what happens when the cone is not convex, in particular for what cones symmetry breaking occurs. 

In the case of constant mean curvature and of overdeterminated problems, symmetry breaking results have been proved in \cite{IPW} for cones spanned by domains $D\subset\S^{N-1}$,  for which the first nontrivial Neumann eigenvalue $\lambda_1(D)$ of the Laplace-Beltrami operator $-\Delta_{\S^{N-1}}$ on $D$ is less than $N-1$. 

One of the main purposes of our paper is to analyze classes of nonlinearities and related classes of spherical sectors $S_D$ for which there exist nonradial positive weak solutions of \eqref{eq:semilinearelliptic}. We will show that the break of symmetry is related to a bound on the Neumann eigenvalue $\lambda_1(D)$. 

A similar result will be established also for a critical Neumann problem in the whole unbounded cone $\Sigma_D$. 

The strategy to get symmetry breaking is to study the Morse index of radial solutions of \eqref{eq:semilinearelliptic} (see Definition \ref{def_Morse}) to show that their instability increases in dependence of the eigenvalue $\lambda_1(D)$. This allows to deduce that, for some nonlinearities,  radial solutions cannot be least-energy solutions of \eqref{eq:semilinearelliptic}, proving so the existence of nonradial positive solutions. Since $S_D$ is a radial domain, to compute the Morse index of a radial solution $\tilde u$ of \eqref{eq:semilinearelliptic}, we decompose the spectrum of the linearized operator $L_{\tilde u}=-\Delta -f'(\tilde u)$ as the sum of  radial eigenvalues of a singular one-dimensional operator and of  Neumann eigenvalues of the Laplace-Beltrami operator $-\Delta_{\S^{N-1}}$ on $D$. In this way we get precise formulas for the Morse index of radial solutions (Proposition \ref{P1}) which are interesting in themselves. Then a crucial role to detect a possible symmetry breaking is played by the first eigenvalue $\hat\Lambda_1^{rad}$ of the following singular eigenvalue problem in the interval $(0,1):$
\beq
\begin{cases}\label{eq:radial}
-\psi''-\frac{N-1}{r}\psi'-f'(\tilde u)\psi=-\frac{\mu}{r^2}\psi &\text{in }(0,1)\\
\psi(1)=0
\end{cases}
\eeq
 where with abuse of notation, $\tilde u(r)=\tilde u(|x|)$.
 
 As explained in Section 2 and Section 3 the number of the negative eigenvalues of \eqref{eq:radial} corresponds to the radial Morse index of $\tilde u$, which we denote by $m_{rad}(\tilde u)$, i.e. to the number of negative eigenvalues of the linearized operator $L_{\tilde u}$ in the space of radial functions $H_{0,rad}^1(S_D\cup\Gamma_1 )\subset H_0^1(S_D\cup \Gamma_1 )$.
 
\begin{remark}
It is important to stress that the eigenvalues of \eqref{eq:radial} do not depend on the domain $D\subset\S^{N-1}$ which spans the cone but only on the nonlinearity $f$ and on the radial solution $\tilde u$ considered. Another way of seeing this, is by observing that a weak radial solution of \eqref{eq:semilinearelliptic} is indeed a radial solution of the corresponding Dirichlet problem in the unit ball $B_1$, i.e. it satisfies:
\beq
\begin{cases}
-\Delta\tilde u=f(\tilde u) & \text{in }B_1\\
\tilde u=0&\text{on }\partial B_1
\end{cases}
\eeq
Therefore the radial Morse index $m^{rad}(\tilde u)$ is the same as the one of $\tilde u$ in the space $H_0^1(B_1)$.
\end{remark}
Our main result about the break of symmetry in the spherical sector $S_D$ is the following.
\begin{theorem}\label{Th:symbreak}
Let us assume that the nonlinearity $f$ is such that:
\begin{enumerate}
\item there exists a unique positive bounded radial solution $\tilde u$ of \eqref{eq:semilinearelliptic},
\item there exists a weak positive solution $\bar u$ of \eqref{eq:semilinearelliptic} with Morse index one, in the space $H_0^1(S_D\cup\Gamma_1 )$.
\end{enumerate}
Then if 
\beq\label{conditions}
\text{either }\hat\Lambda_1^{rad}\ge 0\quad\text{ or }\quad\lambda_1(D)<-\hat \Lambda_1^{rad}
\eeq
then the solution $\bar u$ is not radial.  Here $\lambda_1(D)$ is the first nontrivial Neumann eigenvalue of $-\Delta_{\S^{N-1}}$ on $D$ and $\hat \Lambda_1^{rad}$ is the first eigenvalue of problem \eqref{eq:radial} corresponding to $\tilde u$.
\end{theorem}
Assumptions $1)$ and $2)$ are satisfied by many types of nonlinearities, in particular by the Lane-Emden nonlinearity $f(s)=s^p$, $1<p<\frac{N+2}{N-2}$ (see Section 3).
 Actually, considering $f(s)=s^p$ and studying the behaviour of $\hat\Lambda_1^{rad}$ as $p\to\frac{N+2}{N-2}$ we obtain:
 \begin{theorem}\label{Th:existence}
 Let $S_D$ be a spherical sector spanned by a domain $D\subset \S^{N-1}$ for which $\lambda_1(D)<N-1$. Then there exists $p_0\in \bigg(1,\frac{N+2}{N-2}\bigg)$ such that for every $p\in \bigg(p_0,\frac{N+2}{N-2}\bigg)$ there exists a nonradial positive solution of  \eqref{eq:semilinearelliptic} for $f(u)=u^p$.
 \end{theorem}
Theorem \ref{Th:existence} implies that a radial symmetry result as the one of Gidas-Ni-Nirenberg cannot hold in a sector $S_D$ for which $\lambda_1(D)<N-1$.
 
 It is interesting to observe that the bound $\lambda_1(D)<N-1$ implying the break of symmetry is the same as the one for the constant mean curvature and the overdetermined problems obtained in \cite{IPW}.  Some examples of domains $D\subset\S^{N-1}$ for which $\lambda_1(D)<N-1$ are described in \cite{IPW}.
 
 Finally we also study the question of symmetry of positive solutions to the critical Laplace equation in the unbounded cone $\Sigma_D$. Namely we consider the problem:
 \beq
 \begin{cases}\label{Pbcone}
 -\Delta u=u^{p_s} & \text{ in }\Sigma_D\\
 u>0& \text{ in }\Sigma_D\\
 \frac{\partial u}{\partial\nu}=0&\text{ on }\partial\Sigma_D
 \end{cases}
 \eeq 
where $p_S=\frac{N+2}{N-2}=2^*-1$. 
 
 If the cone is the whole $\R^N$, i.e. $D=\S^{N-1}$, it is well known that all positive solutions of \eqref{Pbcone} are given by the radial function
 \beq\label{Bubble}
 U(x)=\alpha_N\bigg(\frac{1}{1+|x|^2}\bigg)^\frac{N-2}{2}, \quad \alpha_N=(N(N-2))^\frac{N-2}{4},
 \eeq
 as well as by any rescaling or translation of it.
 These functions are usually called standard bubbles. They are also the functions which, up to a constant, achieve the best Sobolev constant for the embedding of the space $D^{1,2}(\R^N)=\{u\in L^{2^*}(\R^N):|\nabla u|\in L^2(\R^N)\}$ into $L^{2^*}(\R^N)$.
 
 In the case of convex cones it has been proved in \cite{LPT} that the standard bubbles are the only positive solutions of \eqref{Pbcone}. This result has been extended to critical equations for more general operators in \cite{CFR}. Thus the question is whether nonradial solutions of \eqref{Pbcone} exist in nonconvex cones.
 
 A first result in this direction has been obtained in \cite{CP} where the existence of a positive nonradial solution of \eqref{Pbcone} is established under some conditions involving the local convexity of $\Sigma_D$ at a boundary point and the measure of $D$ \cite[Corollary 3.5, Theorem 3.6]{CP}.
 
 Here we obtain a more precise characterization of a class of nonconvex cones for which \eqref{Pbcone} admits a nonradial positive solution. Our result is the following:
 \begin{theorem}\label{Th:existencecone}
Let $\Sigma_D$ be a cone spanned by a domain $D\subset \S^{N-1}$ such that
\beq\label{conddomain}
 \bar D\subset\S_{N-1}^+\quad and\quad \lambda_1(D)<N-1
\eeq
 where $\S^+_{N-1}$ is the half-sphere.  Then there exists a positive solution $w$ of \eqref{Pbcone} which is nonradial and fast decaying, i.e. $w(x) = O(|x|^{2-N})$ as $|x|\to +\infty$. 
 \end{theorem}
 Note again the similarity of the bound on $\lambda_1(D)$ in \eqref{conddomain},  Theorem \ref{Th:existence} and the results in \cite{IPW}.
 
 The proof of Theorem \ref{Th:existencecone} relies on a careful analysis of the Morse index of the standard bubble $U$ to show that it becomes more unstable as soon as $\lambda_1(D)$ crosses the value $N-1$ (see Theorem \ref{Th:formula}). This, in turn,  allows to prove that the standard bubbles cannot bee minimizers for the Sobolev quotient $Q_{\Sigma_D}$ defined in \eqref{Sobquot}. On the other side, if $\bar D\subset\S^+_{N-1}$ then a minimizer for $Q_{\Sigma_D}$ should exist by a result of \cite[Theorem 3.3]{CP} so it gives a positive nonradial solution of \eqref{Pbcone}.
 
 The paper is organized as follows.  In Section 2 we describe some preliminary results needed to study the eigenvalues of the linearized operator $L_{\tilde u}=-\Delta-f'(\tilde u)$ at a radial solution $\tilde u$ of \eqref{eq:semilinearelliptic}. In Section 3 we study the Morse index of a radial solution of \eqref{eq:semilinearelliptic} by using a spectral decomposition.  Then we consider some classes of nonlinearities for which least energy positive solutions exist and prove Theorem \ref{Th:symbreak}.
In the same section we analyze the case of Lane-Emden nonlinearities $f(u)=u^p$ obtaining Theorem \ref{Th:existence}. Finally in Section 4 we compute the Morse index of the standard bubble in the unbounded cone $\Sigma_D$ and prove Theorem \ref{Th:existencecone}.

\section{Preliminary results}
To compute the Morse index of a radial solution $\tilde u$ of \eqref{eq:semilinearelliptic} we need to analyze the spectrum of the linearized operator at $\tilde u$ in the space $H_0^1(S_D\cup\Gamma_1 )$. Therefore in this section we consider a general linear operator $L_a$ of the type $-\Delta-a$ and we study its eigenvalues.

Let $a(x)$ be a radial function in $L^\infty(S_D)$ and, for any $v\in H_0^1(S_D\cup\Gamma_1 )$, we consider the linear operator $L_a(v):H_0^1(S_D\cup\Gamma_1 )\to\R$ defined by
\beq\label{def:operator}
L_a(v)\varphi:=\int_{S_D} \left( \nabla v\cdot\nabla\varphi-a(x)v\varphi \right) dx,\qquad \varphi\in H_0^1(S_D\cup\Gamma_1 )
\eeq
and we let $Q_a:H_0^1(S_D\cup\Gamma_1 )\to\R$ be the quadratic form associated to $L_a$, i.e.
\beq\label{def:quadraticform}
Q_a(v):=\int_{S_D} \left(|\nabla v|^2-a(x)v^2 \right) dx.
\eeq

We recall that  the first eigenvalue $\Lambda_1$ of $L_a$, is defined by
\beq\label{def:eigenvalue}
\Lambda_1:=\min_{v\in H_0^1(S_D\cup\Gamma_1 ), v\ne 0} \frac{Q_a(v)}{\int_{S_D}v^2(x)dx},
\eeq
and it is attained at a corresponding eigenfunction $\psi_1\in H_0^1(S_D\cup\Gamma_1 )$ satisfying
\beq\label{def:eigprob}
\begin{cases}
-\Delta\psi_1-a(x)\psi_1=\Lambda_1\psi_1 & \textmd{in } S_D\\
\frac{\partial \psi_1}{\partial \nu}=0 & \textmd{on } \Gamma_1 \\
\psi_1=0& \textmd{on }  D\\
\psi_1>0 &  \textmd{in } S_D.
\end{cases}
\eeq
Then, iteratively for $i \geq 2$, we can define the eigenvalues $\Lambda_i$ by using their min-max characterization 
\beq\label{def:eigenvalues }
\Lambda_i:=\underset{\substack{ v\in H_0^1(S_D\cup\Gamma_1 ),\\ v\ne 0,\\ v\perp\{\psi_1,...,\psi_{i-1}\}}}{\min} \frac{Q_a(v)}{\int_{S_D}v^2(x)dx}=\underset{\substack{W\subset H_0^1(S_D\cup\Gamma_1 ),\\ \dim{W}=i}}{\min} \underset{\substack{ v\in W,\\ v\ne 0,}}{\max} \frac{Q_a(v)}{\int_{S_D}v^2(x)dx} \,,
\eeq
where the condition $v\perp\psi_j$ stands for the orthogonality in $L^2(S_D)$ and $\psi_j$ is a function that attains $\Lambda_j$ for $j=1,...,i-1$. Again the infimum in \eqref{def:eigenvalues } is attained at a function $\psi_i\in H_0^1(S_D\cup\Gamma_1 )$ which is a weak solution to \eqref{def:eigprob}, with $\Lambda_i$ in place of $\Lambda_1$.

%

Before introducing a singular eigenvalue problem associated to $L_a$, it is useful to give a Hardy inequality in a sector. This inequality is probably well-known, but we couldn't find a reference in the literature and we give a proof for the sake of completeness. 

\begin{proposition}[Hardy inequality] \label{prop_Hardy}
Let $N \geq 3$ and let $S_D$ be given by \eqref{def:bddconespanned}. For any $v\in H_0^1(S_D\cup \Gamma_1 )$,  it holds
\begin{equation} \label{hardy}
\frac{(N-2)^2}{4}\int_{S_D}\frac{v^2}{|x|^2}dx\le\int_{S_D}|\nabla v|^2dx.
\end{equation}
\end{proposition}
\begin{proof}
A classical way to prove Hardy type inequalities is to use superharmonic functions. Indeed, if we assume that there exists $G:\overline S_D\to\R$ and $k >0$ such that
$$
\begin{cases}
-\Delta G\ge 0 & \textmd{ in } S_D\\
\partial_\nu G=0 & \textmd{ on } \Gamma_1 \\
G \geq k&\textmd{ in }  S_D
\end{cases}
$$
then an integration by parts implies that
$$
0\le-\int_{S_D} \Delta G\varphi dx= \int_{S_D} \nabla G\cdot\nabla\varphi dx
$$
for any $\varphi\in C_c^1(S_D\cup \Gamma_1 )$ with $\varphi\ge 0$. Given $v \in H_0^1(S_D\cup \Gamma_1 )$, by density we can choose $\varphi=\frac{v^2}{G}$ and we have
$$
0\le \int_{S_D}2v\frac{\nabla G\cdot\nabla v}{G}dx-\int_{S_D}v^2\frac{|\nabla G|^2}{G^2}dx.
$$
Let $\delta >0$, Young inequality yields
$$
\int_{S_D}v^2\frac{|\nabla G|^2}{G^2}dx\le \delta \int_{S_D}v^2\frac{|\nabla G|^2}{G^2}dx+\frac{1}{\delta}\int_{S_D}|\nabla v|^2dx,
$$
that is
\begin{equation} \label{hardy1}
\delta(1-\delta)\int_{S_D}v^2\frac{|\nabla G|^2}{G^2}dx\le\int_{S_D}|\nabla v|^2dx.
\end{equation}
Now we notice that, if $2-N\le \alpha\le 0$ then we can take $G=|x|^\alpha$. Indeed
$$
-\Delta G= \alpha (2-N-\alpha)|x|^{\alpha-2}\ge 0 \quad \text{for}\quad 2-N\le \alpha\le 0
$$
and we also have that $\partial_\nu G = 0 $ on $\Gamma_1 $ and $G \geq 1 $ in $S_D$. Thus from \eqref{hardy1} we obtain
$$
\delta(1-\delta)\alpha^2\int_{S_D}\frac{v^2}{|x|^2}dx\le\int_{S_D}|\nabla v|^2dx.
$$
Since
$$
\delta(1-\delta)\alpha^2\ge \frac{(N-2)^2}{4}
$$
and it is achieved for $\alpha=2-N$ and $\delta=\frac{1}{2}$, then we obtain \eqref{hardy}.
\end{proof}

Now we introduce the following singular eigenvalue problem which will be crucial to prove our main results:
\begin{equation}\label{def:singularproblem}
\begin{cases}
-\Delta \hat\psi-a(x)\hat\psi=\frac{\hat\Lambda}{|x|^2}\hat\psi & \textmd{ in } S_D\\
\partial_\nu \hat\psi =0 & \textmd{ on }  \Gamma_1 \\
\hat\psi=0& \textmd{ on } D.
\end{cases}
\end{equation}
We notice that problem \eqref{def:singularproblem} is well defined in $H_0^1(S_D\cup \Gamma_1 )$ thanks to Hardy inequality \eqref{hardy}, and hence by a weak solution to \eqref{def:singularproblem} we mean $\hat\psi\in H_0^1(S_D\cup\Gamma_1 )$ such that
$$
\int_{S_D}\nabla \hat\psi\cdot \nabla \varphi dx-\int_{S_D}a\hat\psi\varphi dx=\hat\Lambda\int_{S_D}\frac{\hat\psi\varphi}{|x|^2}dx
$$
for every $\varphi\in H_0^1(S_D\cup\Gamma_1 )$.

We start by defining the \emph{singular eigenvalues}, as follows. Let  
\beq \label{hatLambda1}
\hat\Lambda_1:=\inf_{v\in H_0^1(S_D\cup\Gamma_1 ), v\ne 0} \frac{Q_a(v)}{\int_{S_D}|x|^{-2}v^2(x)dx} \,,
\eeq
where $Q_a$ is given by \eqref{def:quadraticform}. Here we stress that the infimum is taken in $H_0^1(S_D\cup\Gamma_1 )$ so that the Hardy inequality assures that the denominator on the right-hand side of \eqref{hatLambda1} is finite. Further properties related to \eqref{hatLambda1} are described in the following proposition.

\begin{proposition} \label{prop_hatLambda1}
Let $\hat\Lambda_1$ be given by \eqref{hatLambda1}, and assume that $\hat\Lambda_1<0$. Then $\hat\Lambda_1$ is attained at a function $\hat\psi_1\in H_0^1(S_D\cup\Gamma_1 )$ which is a positive weak solution to \eqref{def:singularproblem}.  Moreover, if $a\in C^{0,\beta}(S_D)$ for some $0<\beta<1$, then $\hat\psi_1\in C_{loc}^{2,\gamma }(\bar S_D\setminus \{0\})$ for some $0<\gamma <1$ and $\hat\psi_1$ is a classical solution to \eqref{def:singularproblem} in $S_D\setminus\{0\}$ for $\hat\Lambda=\hat{\Lambda_1}$.
\end{proposition}

\begin{proof}
We follow the proof of \cite[Proposition 3.1]{AG}. Let $v_n\in H_0^1(S_D\cup\Gamma_1 )$ be a minimizing sequence for  \eqref{hatLambda1}, and assume that it is such that 
$$
\int_{S_D} v_n^2dx=1.
$$
By definition
\beq\label{eq:minim}
\int_{S_D}(|\nabla v_n|^2-av_n^2)dx=\hat\beta_n\int_{S_D}|x|^{-2}v_n^2dx
\eeq
with $\hat\beta_n\searrow\hat\Lambda_1$ as $n\to\infty$. We first notice that 
\begin{equation} \label{2.13a}
\int_{S_D}|\nabla v_n|^2dx\le \|a\|_{\infty}.
\end{equation}
Indeed, since $\hat\Lambda_i<0$, we can assume that $\hat\beta_n\le 0$ and from \eqref{eq:minim} and \eqref{2.13a} we get
$$
\int_{S_D}|\nabla v_n|^2dx\le \int_{S_D} av_n^2 dx\le \|a\|_{\infty} \,.
$$
Hence, up to a subsequence, $v_n\rightharpoonup \bar v $ weakly in $H_0^1(S_D\cup\Gamma_1 )$ and strongly in $L^2(S_D)$, and in particular
$$
\lim_{n\to\infty}\int_{S_D}av_n^2dx=\int_{S_D}a \bar v^2 dx.
$$
Moreover from \eqref{hardy} and \eqref{2.13a} we have
$$
\int_{S_D}|x|^{-2}v_n^2 dx\le \frac{4}{(n-2)^2}\|a\|_{\infty}  \,.
$$
Now we check that $v$ minimizes the quotient in \eqref{hatLambda1}, namely we prove that
$$
\int_{S_D}\bigg(|\nabla \bar v|^2-a\bar v^2\bigg)dx-\hat\Lambda_1\int_{S_D}|x|^{-2}\bar v^2(x)dx\le 0.
$$
Indeed, since $v_n \to \bar v$ a.e. in $S_D$, Fatou's Lemma yields
$$
\int_{S_D}\bigg(|\nabla \bar v|^2-a\bar v^2\bigg)dx-\hat\Lambda_1\int_{S_D}|x|^{-2}\bar v^2(x)dx\le \liminf_{n\to\infty}\bigg(\int_{S_D}(|\nabla v_n|^2-\beta_n|x|^{-2}v_n^2)dx\bigg)-\int_{S_D}a\bar v^2dx
$$
and \eqref{eq:minim} implies 
$$
\int_{S_D}\bigg(|\nabla \bar v|^2-a\bar v^2\bigg)dx-\hat\Lambda_1\int_{S_D}|x|^{-2}\bar v^2(x)dx\le \liminf_{n\to\infty}\int_{S_D}av_n^2 dx-\int_{S_D}a\bar v^2 dx=0 \,,
$$
where the last inequality follows from the strong convergence $v_n \to \bar v$ in $L^2(S_D)$.

Now, since $\bar v$ minimizes the right-hand side of \eqref{hatLambda1}, we have that $\bar v$ is a weak solution to \eqref{def:singularproblem} for $\hat\Lambda=\hat\Lambda_1$ and then the assertion of the proposition follows by taking $\psi_1=\bar v$. The regularity properties of $\psi_1$ follow from standard elliptic estimates in \cite{GT}.
\end{proof}

Arguing as in Section 3 of \cite{AG} we have that Proposition \ref{prop_hatLambda1} implies that if $\hat\Lambda_1<0$ then we can start an iterative procedure and define the subsequent singular eigenvalues $\hat\Lambda_2,\hat\Lambda_3, \ldots$ at least until they remain negative. Indeed, when $\hat\Lambda_1$ is negative and it is attained at a function $\hat\psi_1\in H_0^1(S_D\cup\Gamma_1 )$, we can define
$$
\hat\Lambda_2:= \inf_{v \in H_0^1(S_D\cup\Gamma_1 ), v\ne 0, v\perp\hat\psi_1}\frac{Q_a(v)}{\int_{S_D}|x|^{-2}v^2(x)dx}.
$$
Iteratively if $\hat\Lambda_{j}<0$ and it is attained at a function $\hat\psi_{j}$, $j=1,\ldots,i-1$, we can define
\beq\label{def:singeigen}
\hat\Lambda_i:=\inf_{v\in H_0^1(S_D\cup\Gamma_1 ), v\ne 0, v\perp\{\hat\psi_1,...,\hat\psi_{i-1}\}}\frac{Q_a(v)}{\int_{S_D}|x|^{-2}v^2(x)dx}.
\eeq
We also notice that if $\hat\psi_i,$ is a function in $H_0^1(S_D\cup \Gamma_1 )$ for which $\hat\Lambda _i$ is attained then it satisfies \eqref{def:singularproblem} for $\hat \Lambda=\hat\Lambda_i$. 

We will say that $\hat\Lambda_i<0$ is a \emph{singular eigenvalue} and the corresponding nontrivial function $\hat\psi_i$ satisfying \eqref{def:singularproblem} will be called \emph{singular eigenfunction}.

Now we prove that the numbers $\hat \Lambda_i$ defined iteratively by \eqref{def:singeigen} are indeed achieved, whenever they are negative. 
\begin{proposition}\label{prop_hatLambdai}
Let $\hat\Lambda_i$ be given by \eqref{def:singeigen} for some $i\in\N$ and assume that $\hat\Lambda_i<0$. Then there exists a function $\hat \psi_i\in H_0^1(S_D\cup\Gamma_1 )$ such that \eqref{def:singeigen} is attained at $\hat \psi_i$, which is a weak solution to \eqref{def:singularproblem} with eigenvalue $\hat\Lambda_i$. Moreover, if $a\in C^{0,\beta}(S_D)$ for some $0<\beta<1$, then $\hat\psi_i\in C_{loc}^{2,\gamma }(\bar S_D\setminus \{0\})$ for some $0<\gamma <1$ and $\hat\psi_i$ is a classical solution to \eqref{def:singularproblem} in $S_D\setminus\{0\}$.
\end{proposition}
\begin{proof}
The proof is analogous to the one of Proposition \ref{prop_hatLambda1}, and for this reason we only give a sketch of the proof enlightening the main differences. Let $v_n$ be a minimizing sequence for \eqref{def:singeigen} such that 
$$
\int_{S_D} v_n^2 dx = 1
$$
and with $v_n\perp\psi_j$, for any $ n \in \mathbb N $ and $j=1,...,i-1$. Hence, $v_n$ converges to a function $\bar v$ weakly in $H_0^1(S_D\cup\Gamma_1 )$, strongly in $L^2(S_D)$ and pointwise a.e. in $S_D$. 

From Hardy inequality \eqref{hardy} and since $\hat \Lambda_i<0$, we have that 
$$
\int_{S_D} |x|^{-2}v_n^2dx\le \frac{4}{(n-2)^2} \|a\|_\infty,
$$
which implies that there exists a subsequence of $v_n$, that we denote $v_{n_k}$,  such that
$$
0=\lim_{k\to\infty} \int_{S_D}|x|^{-2}v_{n_k}\psi_jdx=\int_{S_D}|x|^{-2}\bar v\psi_jdx
$$
for $j=1,...,i-1$, meaning that $\bar v\perp \{\psi_1,...,\psi_{i-1}\}$. Then as before it follows that $\hat \Lambda_i$ is attained and $\bar v=\psi_i$ is a weak solution to \eqref{def:singularproblem} corresponding to $\hat\Lambda_i$. The rest of the proof is completely analogous to the one of Proposition \ref{prop_hatLambda1} and for this reason is omitted.
\end{proof}
%
%
%
%
%
%

Up to now we have introduced two different families of eigenvalues $\Lambda_i$ and $\hat\Lambda_i$, where the last ones are in some sense easier to be studied as we are going to see in the next sections. The crucial point here is that the number of negative eigenvalues $\hat\Lambda_i$ coincides with the number of negative eigenvalues $\Lambda_i$.

\begin{proposition}\label{prop:radial=nonradial}
Let $k_a$ and $\hat k_a$ be defined as:
$$
k_a:=\#\{i\in\N:\Lambda_i<0\}
$$
and
$$
\hat k_a:=\#\{i\in\N :\hat\Lambda_i<0\} \,,
$$
respectively. Then we have that $k_a=\hat k_a$.
\end{proposition}
\begin{proof}
Let $\hat\psi \in H_0^1(S_D\cup\Gamma_1 )$ be an eigenfunction for problem (\ref{def:singularproblem}) corresponding to a singular negative eigenvalue $\hat\Lambda$, i.e.
\begin{equation} \label{hatpsieq}
\begin{cases}
-\Delta\hat\psi-a\hat\psi=\frac{\hat\Lambda}{|x|^2} \hat\psi & \textmd{ in } S_D\\
\frac{\partial \hat\psi}{\partial \nu}=0 & \textmd{ on } \Gamma_1 \\
\hat\psi=0& \textmd{ on } D\\
\end{cases}
\end{equation}
with $\hat \Lambda <0$. Multiplying by $\hat\psi$ and integrating by parts we get
$$
\int_{S_D}(|\nabla\hat\psi|^2-a\hat\psi^2)dx=\hat\Lambda\int_{S_D}\frac{\hat\psi^2}{|x|^2} dx<0
$$
which means that any eigenfunction of \eqref{def:singularproblem} corresponding to a negative eigenvalue makes negative the quadratic form $Q_a$.  Hence 
$$
\hat k_a\le k_a
$$
In order to show the reverse inequality let us suppose by contradiction that
\beq\label{ineq1}
\hat k_a <k_a
\eeq
and denote by $W$ the $k_a-$dimensional space spanned by the eigenfunction of $L_a$ corresponding to the negative eigenvalue $\Lambda_i$.  By definition we get
$$
\hat \Lambda_i\le\max_{v\in W, v\ne 0}\frac{Q_a(v)}{\int_{S_D}|x|^{-2}v^2(x)dx}<0
$$
since any $v\in W$ satisfies $Q_a(v)<0$. This shows that \eqref{ineq1} is not possible.
\end{proof}
Thanks to the geometry of $S_D$ and since the function $a$ is radial we can use separation of variable to study the singular eigenvalue problem \eqref{def:singularproblem}.  Let $H_{0,rad}^1(S_D\cup\Gamma_1 )$ be the subspace of $H_0^1(S_D\cup\Gamma_1 )$ given by radial functions and we define
\beq\label{def:radialeigenvalue1}
\hat \Lambda_1^{rad}:=\inf\bigg\{\frac{Q_a(v)}{\int_{S_D}|x|^{-2}v^2(x)dx}: \ v\in H_{0,rad}^1(S_D\cup\Gamma_1 )\bigg\},
\eeq
and, for $i\ge 2$, 
\beq\label{def:radialeigenvalue}
\hat \Lambda_i^{rad}:=\inf\bigg\{\frac{Q_a(v)}{\int_{S_D}|x|^{-2}v^2(x)dx}:  \ v\in H_{0,rad}^1(S_D\cup\Gamma_1 ) \text{ and } v\perp \psi_1,...,\psi_{i-1}\bigg\}.
\eeq

As for the numbers $\hat\Lambda_i$, previously defined we can prove the following result.
\begin{proposition}
Let $\hat\Lambda_i^{rad}$ be defined by \eqref{def:radprob} and assume that $\hat\Lambda_i^{rad}<0$. Then there exist radial functions $\psi_i^{rad}\in H_0^{1,rad}(S_D\cup\Gamma_1 )$ which achieve \eqref{def:radprob}. Moreover,  using polar coordinates, $\psi_i^{rad}$ are weak solutions of
\beq\label{6'}
\begin{cases}
-(\psi_i^{rad})''-\frac{N-1}{2}\psi_i^{rad '}-a(r)\psi_i^{rad}=\frac{\hat\Lambda_i^{rad}}{r^2}\psi_i^{rad}&\text{for } r\in (0,1)\\
\psi(1)=0
\end{cases}
\eeq
where $a(r)=a(x)$ with $r=|x|$and $'$ denotes the derivative with respect to $r$.
\end{proposition}

Note that the eigenvalue problem \eqref{6'} can be written as the following Sturm Liouville problem 
\beq\label{def:radprob}
\begin{cases}
-(r^{N-1}(\psi_i^{rad})')'-r^{N-1}a(r)\psi_i^{rad}=r^{N-3}\Lambda_i^{rad}\psi_i^{rad} \\
\psi_i^{rad}(1)=0
\end{cases}\quad r\in(0,1).
\eeq
where $\psi_i^{rad}\in H_0^{1,rad}(S_D\cup\Gamma_1 )$.

Regarding the angular component on $D$, we denote by $\Delta_{\S^{N-1}}$ the Neumann Laplace-Beltrami operator on $\S^{N-1}$ and consider the Neumann eigenvalues of $-\Delta_{\S^{N-1}}$ on the domain $D$. It is well-known that $(-\Delta_{\S^{N-1}})^{-1}$ is compact and selfadjoint in $L^2(D)$ and admits a sequence of eigenvalues 
\begin{equation} \label{lambda}
0=\lambda_1<\lambda_2\le \ldots \lambda_j \leq \ldots
\end{equation} 
and corresponding eigenfunctions $Y_j(\theta) \in L^2(D)$ (where $\theta$ is the system of coordinates on $D$ induced by the spherical coordinates in $\R^N$) which form a Hilbert basis for $L^2(D)$ and such that
\beq\label{def:lapleig}
-\Delta_{\S^{N-1}}Y_j(\theta)=\lambda_j Y_j(\theta)\quad \theta\in D 
\eeq
and 
$$
\int_{S_D} \nabla_\theta Y_j(\theta) \cdot \nabla_\theta Y_i(\theta) d\sigma(\theta) = 0 \quad \text{ for } i \neq j \,,
$$
where we denote by $\nabla_\theta$ the gradient with respect to $\theta \in D$.

Since $a$ is a radial function, the singular eigenvalue $\hat \Lambda_i$ can be decomposed as follows.
\begin{proposition}\label{prop:decomp}
Let  $\hat \Lambda_i$ and $\lambda_j$ be given by \eqref{def:singeigen} and \eqref{lambda}, respectively. If $\hat \Lambda_i<0$ for some $i \in \mathbb{N}$ then there exist $ k\ge 1 $ and $j\ge 0$ such that 
\begin{equation}\label{def:decompeigenv}
\hat\Lambda_i=\hat\Lambda_k^{rad}+\lambda_j \,.
\end{equation}
Viceversa, let us assume that there exist $k \geq 1$ such that $\hat\Lambda_k^{rad}<-\lambda_j$ for some $j\ge 0$,  then $\hat\Lambda_i$ given by \eqref{def:decompeigenv} is a negative singular eigenvalue for \eqref{def:singularproblem}.
\end{proposition}
\begin{proof}
Let $\hat\Lambda<0$ and let $\hat \psi\in H_0^1(S_D\cup\Gamma_1 )$ be an eigenfunction associated to $\hat\Lambda$, i.e. $\hat \psi $ is a solution to \eqref{hatpsieq}. Since $\{Y_j\}$ are an orthonormal base on $L^2(D)$, we can write
$$
\hat\psi(r,\theta)=\sum_{j=0}^\infty \hat\psi_j(r)Y_j(\theta)
$$
for $r\in(0,1]$ and $\theta\in D\subset \S^{N-1}$, where
\beq\label{eq:eigenf}
\hat\psi_j(r):=\int_{D}\hat\psi(r,\theta)Y_j(\theta)d\sigma(\theta) \,.
\eeq
We first notice that, since $\hat \psi\in H_0^1(S_D\cup\Gamma_1 )$, from Parseval identity and Hardy inequality we have that
\begin{equation} \label{hatpsi_jfinite}
\frac{(n-2)^2}{4} \sum_{j=0}^\infty \int_0^1 r^{N-3}\hat\psi_j^2dr \leq \sum_{j=0}^\infty \int_0^1 r^{N-1}(\hat \psi_j')^2dr <+\infty.
\end{equation}
Since $\hat \psi \not\equiv 0$, there exists $j \geq 0$ such that $\hat\psi_j(r)\not \equiv 0$ and we can write
$$
\int_0^1 r^{N-1}\hat\psi_j'\varphi'dr=\int_0^1\int_D r^{N-1} \hat\psi' Y_j(\theta)\varphi' drd\sigma(\theta)=\int_0^1\int_D r^{N-1} \hat\psi'\bigg(Y_j(\theta)\varphi\bigg)' drd\sigma(\theta)
$$
for every $\varphi\in H_0^1(S_D\cup\Gamma_1 )$. We recall that $\nabla \hat \psi=(\hat\psi', r^{-2} \nabla_\theta \hat\psi)$; since $\hat \psi$ is a solution to $\eqref{def:singularproblem}$ and $a$ is radial, we have
\begin{multline*}
\int_0^1 r^{N-1}\hat\psi_j'\varphi'dr=-\int_0^1\int_D r^{N-3}\nabla_\theta\hat\psi\cdot \nabla_\theta (Y_j(\theta)\varphi) drd\sigma(\theta) \\ +\int_0^1\int_D r^{N-1}a\hat\psi Y_j(\theta)\varphi drd\sigma(\theta)  +\hat\Lambda\int_0^1\int_D r^{N-3}\hat \psi Y_j(\theta)\varphi drd\sigma(\theta) \,.
\end{multline*}
From \eqref{eq:eigenf}, \eqref{def:lapleig} and by using integration by parts, we obtain
\begin{flalign*}
\int_0^1 r^{N-1}\hat \psi_j'\varphi'dr=&-\int_0^1r^{N-3}dr \int_D \nabla_\theta(\hat \psi\varphi)\cdot \nabla_\theta Y_j(\theta) d\sigma(\theta)+\int_0^1 r^{N-1}a(r)\hat \psi_j\varphi dr+\hat\Lambda\int_0^1 r^{N-3}\hat \psi_j\varphi dr\\
=&-\lambda_j\int_0^1r^{N-3}\hat \psi_j\varphi dr+\int_0^1 r^{N-1}a(r)\hat \psi_j\varphi dr+\hat\Lambda\int_0^1 r^{N-3}\hat \psi_j\varphi dr\\
=&\int_0^1 r^{N-1}a(r)\hat \psi_j\varphi dr+(\hat\Lambda-\lambda_j)\int_0^1 r^{N-3}\hat \psi_j\varphi dr \,.
\end{flalign*}
Hence we have proved that $\hat \psi_j$ is a weak solution to 
\beq\label{def:probSturmLiouv}
\begin{cases}
-(r^{n-1}(\hat\psi_j)')'-r^{n-1}a\hat\psi_j=r^{N-3}\hat\Lambda_j^{rad}\hat\psi_j & \text{ for } r\in(0,1)\\
\hat\psi_j\in H_0^{1,rad}(D)
\end{cases}
\eeq
i.e.  $\hat\psi_j$ is an eigenfunction of \eqref{def:singularproblem} corresponding to the eigenvalue $(\hat\Lambda-\lambda_j)$.

Since  $\hat\Lambda-\lambda_j<0$,  from \eqref{def:probSturmLiouv} and \cite[Proposition 3.10]{AG} we have that $\hat\Lambda-\lambda_j$ is a radial singular eigenvalue for $L_a$ given by \eqref{def:operator}.

The reverse implication holds as well, namely if $\hat\Lambda_k^{rad}+\lambda_j<0$ for some radial singular eigenvalue $\hat\Lambda_k^{rad}$ with associated eigenfunction $\hat \psi_k^{rad}\in H_0^{1,rad}(S_D\cup\Gamma_1 )$ and for some $\lambda_j$, then the function 
$$
\Psi:=\hat \psi_k^{rad}(r)Y_j(\theta)
$$ 
is such that $\Psi \in H_0^1(S_D\cup\Gamma_1 )$.  
Moreover $\Psi$ weakly solves \eqref{def:singularproblem} corresponding to $\hat\Lambda=\hat\Lambda_k^{rad}+\lambda_j<0$. Indeed for any $\varphi\in H_0^1(S_D\cup\Gamma_1 )$ we have
\begin{flalign*}
\int_{S_D}\nabla\Psi\cdot\nabla\varphi dx=&\int_{S_D}\Psi'\varphi' +\frac{1}{r^2}\nabla_\theta\Psi\cdot\nabla_\theta\varphi dx=\\
=&\int_D Y_j(\theta)d\sigma(\theta)\int_0^1r^{N-1}(\psi^{rad}_k)' \varphi'dr+\int_0^1 r^{N-3}\psi_k^{rad}dr\int_D\nabla_\theta Y_j \cdot\nabla_\theta\varphi d\sigma(\theta)
\end{flalign*}
and from \eqref{def:probSturmLiouv} and \eqref{def:lapleig} it follows
\begin{flalign*}
\int_{S_D}\nabla\Psi\cdot\nabla\varphi dx=&\int_DY_j(\theta)d\sigma(\theta)\int_0^1r^{N-1}(a+\frac{\hat\Lambda_k^{rad}}{r^2})\psi_k^{rad}\varphi dr+\int_0^1 r^{N-3}\psi_k^{rad}dr\lambda_j\int_D Y_j(\theta )\varphi d\sigma(\theta)\\
=&\int_{S_D}a\Psi\varphi+\frac{\hat\Lambda_k^{rad}+\lambda_j}{|x|^2}\Psi\varphi \,dx,
\end{flalign*}
and then $\Psi$ weakly solves \eqref{def:singularproblem} corresponding to $\hat\Lambda=\hat\Lambda_k^{rad}+\lambda_j<0$, which completes the proof.
\end{proof}
\begin{remark}
The results of this section, obtained for negative eigenvalues, also hold for positive eigenvalues which are smaller than $(N-2)^2/4$ (see \cite{AG}). Since we are interested in negative eigenvalues to study the Morse index of solutions we have preferred to avoid further technicalities.
\end{remark}

\section{Application to semilinear elliptic problems}
In this section we study the semilinear elliptic problem \eqref{eq:semilinearelliptic} to the aim of proving the symmetry breaking result in Theorems \ref{Th:existence} and \ref{Th:existencecone}.  We start by showing that weak solutions of \eqref{eq:semilinearelliptic} are bounded.

The following theorem is a classical result on interior boundedness of solutions to elliptic PDEs, which was proved by Serrin in \cite{Serrin}. This result can be extended up to the boundary, whenever the boundary fulfills some suitable regularity. We will need this result at some points in the paper, and in particular it is needed to show boundedness of solutions in a neighborhood of the vertex. Since we were not able to find a reference for this case, we prefer to give a sketch of the proof below.

\begin{proposition}\label{prop:bdd}
Let $u\in H_0^1(S_D\cup\Gamma)$ be a weak solution of
\begin{equation*}
\begin{cases}
-\Delta u=g(x,u)  &\text{in }S_D\\
\partial_\nu u=0  &\text{on } \Gamma \\
u=0  &\text{on } D \\
u>0 &\text{in }S_D
\end{cases}
\end{equation*}
where $g:S_D \times \mathbb R \to \mathbb R$ is such that 
\begin{equation} \label{g_cond}
|g(x,u)| \leq  c_1(x) |u| + c_2(x) \,,
\end{equation}
with $c_1,c_2 \in L^{\frac{N}{2-\epsilon}}(S_D)$, for some $\epsilon>0$ small enough. Then $u$ is bounded in $S_D$ and we have
\beq\label{num}
\|u\|_{L^\infty( S_D)} \leq C \{\|u\|_{L^{2}(S_D)}+\|c_2\|_{L^\infty(S_D)}\} \,,
\eeq
with $C=C(N,\epsilon, M)$, where 
$
M=\|c_1\|_{L^{\frac{N}{2-\epsilon}}(S_D)} \,.
$

In particular, if $u\in H_0^1(S_D\cup\Gamma)$ is a weak solution of 
\begin{equation} \label{pb_app_serrin}
\begin{cases}
-\Delta u=f(u)  &\text{in }S_D\\
\partial_\nu u=0  &\text{on } \Gamma \\
u=0  &\text{on } D \\
u>0 &\text{in }S_D
\end{cases}
\end{equation}
with $|f(s)|\le c|s|^p+b(x)$, $p<\frac{N+2}{N-2}$,  $b(x)\in L^{\frac{N}{2-\epsilon}}$,$ N\ge 3$ then $u$ is bounded in $\overline S_D$.
\end{proposition}
\begin{proof}
As already mentioned, this proof strictly follows the lines of \cite[Theorem 1]{Serrin}. Here, we just want to emphasize that, thanks to the Neumann boundary on $\Gamma$, all the argument can be easily adapted. 

We notice that the boundedness of the solution at interior points follows directly from \cite[Theorem 1]{Serrin}. Boundedness at points on $\partial S_D \setminus \{O\}$ can be obtained by applying standard reflection methods and again using \cite[Theorem 1]{Serrin}. For this reason in the following we prove the assertion only in a neighborhood of the vertex $O$. 

In order to lighten the notation, we set
$$
\|u\|_{p,R}:=\|u\|_{L^p(\Sigma_D \cap B_R)} 
$$
for $p\in[1,\infty]$.

Let $0<R<1/4$ and set
$$\bar u=u+\|c_2\|_{\infty,2R}.
$$
For fixed numbers $q\ge 1$ and $l>\|c_2\|_{\infty,2R}$, we define $F(\bar u)\in C^1(\R)$
$$ F(\bar u)=\begin{cases}
\bar u^q& \textmd{ for }  \|c_2\|_{\infty,2R}\le\bar u\le l\\
ql^{q-1}\bar u-(q-1)l^q& \textmd{ for } l\le \bar u,
\end{cases}
$$
and
$$
G(\bar u)=F(\bar u)F'(\bar u)-q\|c_2\|_{\infty,2R}^{2q-1}.
$$
Notice that 
$$
G'=\begin{cases}
q^{-1}(2q-1)(F')^2 &\bar u<l-\|c_2\|_{\infty,2R}\\
(F')^2&\bar u>l-\|c_2\|_{\infty,2R} \,.
\end{cases}
$$
Let $\eta \in C_c^\infty(B_{2R}(O))$ be such that $\eta \geq 0$, and let
$$\phi=\eta^2G(\bar u).
$$
We notice that we have
$$
\nabla \phi=2\eta\nabla\eta G+\eta^2G'\nabla \bar u \quad  \textmd{ at points where } \{\bar u\ne l-\|c_2\|_{\infty,2R}\}
$$
and then, by using $|G|\le FF'$,  one has
\begin{equation*}
\begin{split}
\nabla\phi\cdot\nabla \bar u-\phi |g| & = \eta^2 G' |\nabla \bar u|^2+2\eta G\nabla \eta\cdot \nabla \bar u-\eta^2G|g| \\
& \ge  \eta^2 (F')^2 |\nabla \bar u|^2-2\eta FF'|\nabla \eta|| \nabla \bar u|-(|c_1|+1)\bar u\eta^2FF' \,.
\end{split}
\end{equation*}
Let 
$$
v:=F(\bar u) \,.
$$ 
Since $\bar uF'=qF$ then we can write 
\beq\label{1}
\nabla\phi\cdot\nabla \bar u-\phi |g| \ge  \eta^2  |\nabla v|^2-2\eta v |\nabla \eta|| \nabla v|-q(|c_1|+1)(\eta v)^2.
\eeq
In the set where $\bar u=l$ we have $\nabla\phi=2\eta\nabla\eta G$ and $\nabla \bar u\equiv 0$, a.e.  so that \eqref{1} holds also on this set a.e.. We may integrate \eqref{1} over $\Sigma_D\cap B_{2R}$ and get
\begin{equation} \label{S1}
\|\eta\nabla v\|_{2,2R}^2\le 2\int_{\Sigma_D\cap B_{2R}} \eta v |\nabla \eta|| \nabla v|\, dx+q\int_{\Sigma_D\cap B_{2R}} (|c_1|+1)(\eta v)^2\,dx.
\end{equation}
We apply H\"older inequality on the r.h.s. and we get
\begin{equation} \label{S2}
\int_{\Sigma_D\cap B_{2R}} \eta v |\nabla \eta|| \nabla v|\, dx\le \|v\nabla \eta \|_{2,2R}\|\eta \nabla v\|_{2,2R}
\end{equation}
$$
\int_{\Sigma_D\cap B_{2R}}( |c_1|+1)(\eta v)^2\, dx =\int_{\Sigma_D\cap B_{2R}} (|c_1|+1)(\eta v)^\epsilon(\eta v)^{2-\epsilon}\, dx\le  \||c_1|+1\|_{N/(2-\epsilon),2R}\| \eta v\|_{2,2R}^\epsilon\| \eta v\|_{2^*,2R}^{2-\epsilon}.$$
Then,  since $c_1\in L^{\frac{N}{2-\epsilon}}(S_D)$ and $S_D$ is bounded, we have: 
$$
 \||c_1|+1\|_{N/(2-\epsilon),2R}\le C(M,N,\epsilon),
$$
which implies
\begin{equation} \label{S3}
\int_{\Sigma_D\cap B_{2R}} (|c_1|+1)(\eta v)^2\, dx\le C\| \eta v\|_{2,2R}^\epsilon \{\|\nabla\eta v\|_{2,2R}^{2-\epsilon}+\| \eta \nabla v\|_{2,2R}^{2-\epsilon}\}
\end{equation}
and from \eqref{S1}-\eqref{S3} we get
\begin{equation} \label{S4}
\|\eta\nabla v\|_{2,2R}^2\le 2\|v\nabla \eta \|_{2,2R}\|\eta \nabla v\|_{2,2R}+qC\| \eta v\|_{2,2R}^\epsilon \{\|\nabla\eta v\|_{2,2R}^{2-\epsilon}+\| \eta \nabla v\|_{2,2R}^{2-\epsilon}\} \,.
\end{equation}
We divide by $ \|v\nabla \eta \|_{2,2R}^2$ and call
$$
z=\|\eta\nabla v\|_{2,2R}/\|v\nabla \eta \|_{2,2R}\qquad\xi= \| \eta  v\|_{2,2R}/\|v\nabla \eta \|_{2,2R} \,;
$$
hence \eqref{S4} can be written as
$$
z^2\le 2z+qC\xi^\epsilon+qC\xi^\epsilon z^{2-\epsilon}.
$$
Lemma 2 in \cite{Serrin} yields
$$
z\le C[2+(Cq\xi^\epsilon)^{1/2}+(Cq\xi^\epsilon)^{1/\epsilon}]\le Cq^{2/\epsilon}(1+\xi),
$$
i.e. 
\beq\label{2}
\|\eta \nabla v\|_{2,2R}\le Cq^{2/\epsilon}(\|\eta v\|_{2,2R}+\|\nabla \eta v\|_{2,2R}).
\eeq
By using Sobolev inequality (see \cite{Adams}) we obtain
\beq\label{3}
\|\eta v\|_{2^*,2R}\le C q^{2/\epsilon}(\|\eta v\|_{2,2R}+\|\nabla \eta v\|_{2,2R}).
\eeq
Equations \eqref{2} and \eqref{3} permit to start a Moser iteration and get the assertion. Indeed, let $h,h' \in \mathbb R$ be such that $h<h'\le 2R$. We choose $\eta$ such that $\eta=1$ in $B_{h'}$, $0\le\eta\le 1$ in $B_h$, $\eta$ identically zero outside $B_h$, and such that $|\nabla \eta| \leq 2(h-h')^{-1}$. By using $\eta$ in \eqref{2} and \eqref{3} yields 
\beq\label{4}
\| \nabla v\|_{2,h'}\le Cq^{2/\epsilon}(h-h')^{-1}\| v\|_{2,h} \,,
\eeq
\beq\label{5}
\|v\|_{2^*,h'}\le Cq^{2/\epsilon}(h-h')^{-1}\| v\|_{2,h} \,.
\eeq
Let $l\to\infty$ in \eqref{5}. Since $v\to \bar u^q$, Lebesgue's monotone convergence theorem yields
\beq\label{6}
\| \bar u^q\|_{2^*,h'}\le Cq^{2/\epsilon}(h-h')^{-1}\| \bar u^q\|_{2,h} \,,
\eeq
i.e.
\begin{equation}\label{7}
\| \bar u\|_{2^*q,h'}\le[ Cq^{2/\epsilon}(h-h')^{-1}]^{1/q}\| \bar u\|_{2q,h}.
\end{equation}
Now we call 
$$
t:=2^*/2\qquad p_\nu:=2t^\nu\quad \nu=0,1,2,...
$$
and $h_\nu=R(1+2^{-\nu})$, $h'_\nu=h_{\nu+1}$ whence \eqref{7} becomes
\begin{equation*}
\| \bar u\|_{p_{\nu+1,h_{\nu+1}}}\le[ C(2t)^{2\nu/\epsilon}]^{1/x^\nu}\| \bar u\|_{p_\nu,h_\nu}.
\end{equation*}
Iteration yields
\begin{equation*}
 \| \bar u\|_{p_{\nu+1,h_{\nu+1}}}\le C^{\sum\limits_{t=1}^\nu 1/t^\nu}(2t)^{(2/\epsilon)(\sum\limits_{t=1}^\nu\nu/t^\nu)}\| \bar u\|_{2,2R}\le C\|\bar u\|_{2,2R},
\end{equation*}
since both series appearing in the last display are convergent. Letting $\nu\to\infty$ and observing that $\|\bar u\|_{\infty,R}\le \lim \|\bar u\|_{p_\nu,h_\nu}$, we find
$$
\|\bar u\|_{\infty, R}\le C\|\bar u\|_{2,2R},
$$
and from
$$
\bar u=u+\|c_2\|_{\infty,2R}
$$
we obtain
$$
\|u\|_{\infty, R}\le C\{\|u\|_{2,2R}+\|c_2\|_{\infty,2R}\}.
$$
This proves \eqref{num} and the first part of the theorem.

In order to prove the latter part of the assertion, we notice that if $u$ is a solution to \eqref{pb_app_serrin} then $u$ satisfies
$$
-\Delta u =c(x)(1+ u) \,, 
$$
with
$$
c(x)= \frac{f(u(x))}{1+u}
$$
with $c(x)$ satisfying
$$
 |c(x)|\le|f(u(x))|\le c|u|^{p-1}|u|+|b(x)|\,.
$$
Since $u \in W^{1,2}$ then $u \in L^{2^*}$ and if we consider $\epsilon<\frac{N-2}{2}(p_s-p)$ it holds
$$
(p-1)\frac{N}{2-\epsilon}<2^*
$$
namely $c |u|^{p-1} \in L^{\frac{N}{2-\epsilon}} $.  Since $b(x)\in L^{\frac{N}{2-\epsilon}}$ then we can apply the first part of the theorem and conclude.

\end{proof}

\begin{remark}\label{lastrmk}
Theorem \ref{prop:bdd} still holds when $g(x,u)\le c (1+|u|^p)$ with $1<p\le p_S$. The critical case $p=p_S$ requires a modification of the proof which can be found in \cite[Appendix D]{P}. 

The extension to a convex cone can be found in \cite[Lemma 2.1]{CFR}. We mention that the convexity of the cone is not needed in the proof. 
\end{remark}

Since the nonlinearity $f$ in \eqref{eq:semilinearelliptic} is locally Lipschitz continuous,  and if we have that the solution $u$ is bounded, then we have that the function $a(x)=f'(u(x))$ belong to $  L^\infty (S_D)$ and we can apply the results of Section 2 to the linearized operator at $u$:
\begin{equation} \label{Lu_def}
L_u(v):=-\Delta v -f'(u) v, \quad v\in  H_0^1(S_D\cup \Gamma_1 )
\end{equation}
defined as in \eqref{def:operator}.
\subsection{Morse index of radial solutions}
To define the Morse index of a solution $u$ to \eqref{eq:semilinearelliptic} we consider the quadratic form associated to $L_u$:
\begin{equation} \label{Qu_def}
Q_u(v):=\int_{S_D}\bigg(|\nabla v|^2-f'(u)v^2\bigg)dx \,, \quad v\in  H_0^1(S_D\cup \Gamma_1 )
\end{equation}
\begin{definition} \label{def_Morse} 
 Let  $u\in H_0^1(S_D\cup \Gamma_1 )$ be a solution of \eqref{eq:semilinearelliptic}. We say that:
 \begin{enumerate}[label=\roman*)]
 \item $u$ is stable (or $u$ has zero Morse index) if $Q_u(w)\ge 0$ for any $w\in C_c^1(S_D\cup\Gamma_1 )$;
 \item $u$ has Morse index equal to the integer $m(u)\ge 1$ if $m(u)$ is the maximal dimension of a subspace of $C^1_c(S_D\cup\Gamma_1 )$ where the quadratic form $Q_u$ is negative definite;
 \item $u$ has infinite Morse index if, for any integer $k\ge 1$, there exists a k-dimensional subspace of $C^1_c(S_D\cup\Gamma_1 )$ where $Q_u$ is negative definite.
 \end{enumerate}
 \end{definition}

Since $L_u$ is a linear compact operator,  the maximal dimension of a subspace of $H_0^1(S_D\cup \Gamma_1 )$ in which $Q_u$
is negative defined is equivalent to the number of negative eigenvalues of $L_u$ in $H_0^1(S_D\cup \Gamma_1 )$, counted with their multiplicity.

From Proposition \ref{prop:radial=nonradial} it follows that the number of negative eigenvalue, counted with their multiplicity,  is equal to the number of negative singular eigenvalues $\hat\Lambda_i$ of the associated singular problem. 
Hence, if $\tilde u\in H_0^{1,rad}(S_D\cup\Gamma_1 )$ is a solution to \eqref{eq:semilinearelliptic}, then we have 
$$
m(\tilde u)=\#\{i\ge 1|\hat\Lambda_i<0\} \,.
$$
We remember that Proposition \ref{prop:decomp} yields
$$
\hat\Lambda_i=\hat\Lambda_k^{rad}+\lambda_j, 
$$
for some $k\ge 1, j\ge 0$,  so that
\beq \label{def:firstMorseindex}
m( \tilde u)=\#\{j\ge 0,k\ge 1|\hat\Lambda_k^{rad}+\lambda_j<0\}.
\eeq
For a radial solution $\tilde u$ of \eqref{eq:semilinearelliptic},  we denote by $m^{rad}(\tilde u)$ the Morse index of $\tilde u$ in the space $H_{0,rad}^1(S_D\cup\Gamma_1)$, i.e.  $m^{rad}(\tilde u)$ is the number of the negative eigenvalues of $L_u$ to which there corresponds a radial eigenfunction.

To show the break of symmetry results it is important to understand how Morse index of $\tilde u$ changes passing from the space $H_{0,rad}^1(S_D\cup\Gamma_1)$ to the space $H_{0}^1(S_D\cup\Gamma_1)$. Next theorem provides precise formulas to get $m(\tilde u)$, depending on the eigenvalues $\lambda_j$ of the Laplace Beltrami operator $-\Delta_{\S^{N-1}}$ on $D$.

\begin{proposition}\label{P1}
Let $\tilde u$ be a bounded weak radial solution of \eqref{eq:semilinearelliptic} and $\hat\Lambda_k^{rad}$ the singular eigenvalue defined in \eqref{def:radialeigenvalue1} and \eqref{def:radialeigenvalue} for the linearized operator  $L_{\tilde u}$, i.e.  considering $a(x)=f'(\tilde u(x))$.  Then
\begin{enumerate}[label=\roman*)]
\item  if $m^{rad}(\tilde u)=0$ then $m(\tilde u)=0$,
\item  if $m^{rad}(\tilde u)=1$ then
\beq\label{morseindex1}
m(\tilde u)= \# 	\{j\ge 1:\lambda_j< -\hat\Lambda^{rad}_1\}+1.
\eeq
\end{enumerate}
More generally it holds 
\begin{enumerate}[label=\roman*), start=3]
\item  if $m^{rad}(\tilde u)=d$ then
$$
m( \tilde u)= \sum_{k=1}^d(d-k+1) \# \{j\ge 1: -\hat\Lambda^{rad}_{k+1}\le\lambda_j< -\hat\Lambda^{rad}_{k}\}+d.
$$
\end{enumerate}
\end{proposition}
\begin{proof}\hfill\\
\begin{enumerate}[label=\roman*)]
\item The proof is an immediate consequence of Proposition \ref{prop:decomp}. Indeed, since $\lambda_j\ge 0$ for any $j\ge 0$, then $\hat\Lambda_k^{rad}\ge 0$ for any $k\ge 1$ and we have
$$
\hat\Lambda_i=\hat\Lambda_k^{rad}+\lambda_j\ge 0
$$
for any $i\ge 1$, which implies that the Morse index of $\tilde u$ is zero.
\item Since $ \tilde u$ is a solution with $m^{rad}( \tilde u)=1$ we have
$$
\hat\Lambda_1^{rad}<0\quad \text{ and }\hat\Lambda_i^{rad}\ge 0\quad \forall i\ge2 \,.
$$
Since $\lambda_0=0$ then we have
$$
m( \tilde u)=\#\{i\ge 1:\hat\Lambda_i< 0\}= \# 	\{j\ge 0:\lambda_j< -\hat\Lambda^{rad}_1\}= \# 	\{j\ge 1:\lambda_j< -\hat\Lambda^{rad}_1\}+1\ge 1 \,,
$$
which is \eqref{morseindex1}.
\item We notice that
$$
m( \tilde u)=\#\{i\ge 1:\hat\Lambda_i< 0\}= \# 	\{k\ge 1,j\ge 0:\lambda_j< -\hat\Lambda^{rad}_k\}=\sum_{k=1}^d \# 	\{j\ge 0:\lambda_j< -\hat\Lambda^{rad}_k\} \,.
$$
Since $\lambda_0=0$ we have
$$
\sum_{k=1}^d \# 	\{j\ge 0:\lambda_j< -\hat\Lambda^{rad}_k\}=\sum_{k=1}^d \# 	\{j\ge 1:\lambda_j< -\hat\Lambda^{rad}_k\}+d \,,
$$
and from $\lambda_0\le \lambda_1\le \lambda_2\le\cdots$, we obtain
$$
\sum_{k=1}^d \# 	\{j\ge 1:\lambda_j< -\hat\Lambda^{rad}_k\}+d=\sum_{k=1}^d(d-k+1) \# \{j\ge 1: -\hat\Lambda^{rad}_{k+1}\le\lambda_j< -\hat\Lambda^{rad}_{k}\}+d,
$$
which ends the proof.
\end{enumerate}
\end{proof}

In the following proposition we give an estimate of the first singular radial eigenvalue $\hat\Lambda_1^{rad}$.

\begin{proposition} \label{autres}
If $\tilde u$ is a bounded positive radial weak solution of \eqref{eq:semilinearelliptic}, then 
$$
\hat\Lambda_1^{rad}>-(N-1).
$$
\end{proposition}

\begin{proof}
Let $\tilde u\in H_0^{1,rad}(S_D\cup\Gamma_1 )$ be a radial solution to \eqref{eq:semilinearelliptic},  with an abuse of notation we write $\tilde u(r)=\tilde u(|x|)$ and define $\eta(r):=\partial \tilde u'(r)$. Notice that $\eta (r) <0$ for $r>0$. Since
$$
-\tilde u''-\frac{N-1}{r}\tilde u'=f(\tilde u)
$$
then it follows that $\eta$ satisfies
$$
-\eta''+\frac{N-1}{r^2}\eta-\frac{N-1}{r}\eta'=f'(\tilde u)\eta \,,
$$
i.e. $\eta(r)<0$ for $r \in(0,1)$ and 
$$
\begin{cases}
(r^{n-1}\eta')'+r^{n-1}f'(\tilde u)\eta=    r^{n-3}(N-1)\eta & \textmd{for } r \in (0,1) \\
\eta(0)=0 \,.& 
\end{cases}
$$
Let $v$ be an eigenfunction associated to the eigenvalue $\hat\Lambda_1^{rad}$ for the following problem
$$
\begin{cases}
(r^{n-1}(v)')' + r^{n-1}f'(\tilde u)v=-r^{n-3}\hat\Lambda_1^{rad}v& \textmd{for } r\in(0,1)\\
v(0)=v(1)=0  \\
v>0
\end{cases}
$$
Assume that $\hat\Lambda_1^{rad}\le -(N-1)$. If $\hat\Lambda_1^{rad}=-(N-1)$, then $\eta$ and $v$ are two solutions of the same Sturm-Liouville problem and they are linearly independent, because $\eta(1)\ne 0=v(1)$. As a consequence of the Sturm separation theorem \cite{Sturm}, the zeros of $\eta$ and $v$ must alternate.  But this leads to a contradiction, since $\eta$ should be zero at a point $\bar r\in(0,1)$, which contradicts $\eta<0$.  If $\hat\Lambda_1^{rad}<-(N-1)$, then by the Sturm-Picone comparison theorem \cite{Picone}, $\eta$ must have a zero between any two consecutive zeroes of $v$.  This leads again to the same contradiction as before. 
\end{proof}

\subsection{Solutions of Morse index one and proof of Theorem \ref{Th:symbreak} }
Before proving Theorem \ref{Th:symbreak} which concerns solutions of Morse index one let us consider a class of nonlinearities for which such solution exists. 
We assume that $f=f(s)$ satisfies the following conditions:
\begin{enumerate}\label{hypNehari}
\item[(F1)] $f(s):\R\to\R$ belongs to $C^1(\R)$;
\item[(F2)] there exists $a_1\in L^\frac{2N}{N+2}(S_D)$ and $a_2>0$ such that
$$
|f(s)|\le a_1(x)+a_2|s|^p\quad\forall s\in \R
$$
for $1<p<p_S $, $p_S=\frac{N+2}{N-2}$ if $N\ge 3$;
\item[(F3)] $f(s)=o(|s|)$ as $s\to 0$;
\item[(F4)] $\exists\alpha>2, r\ge 0$ such that for $|s|\ge r$
$$
0<\alpha F(s)\le sf(s)
$$
where $F$ is the primitive of $f$;
\item[(F5)] $$
\frac{\partial f}{\partial s}(s)>\frac{f(s)}{s}\quad \forall s\in \R\setminus\{0\} \,.
$$
\end{enumerate}
We denote by $J_{S_D}$ the functional
$$
J_{S_D}[u] =\frac{1}{2}\int_{S_D}|\nabla u|^2dx-\frac{1}{2^*}\int_{S_D} |u|^{2^*}dx,\quad u\in H_0^1(S_D\cup\Gamma_1 )
$$
and define the associated Nehari manifold $\mathcal{N}$ as follow
$$
\mathcal{N}:=\bigg\{v\in H_0^{1}(S_D\cup\Gamma_1 ):\langle J'(v),v\rangle=\int_{S_D}|\nabla v|^2dx-\int_{S_D}f(v)v\,dx=0\bigg\}.
$$
We have the following result \cite[Proposition 3.5]{DP}.
\begin{proposition}
Let $f$ satisfy $(F1)-(F5)$. Then $\mathcal{N}$ is a $C^1-$Hilbert manifold of codimension one and:
\begin{enumerate}
\item[(i)] there exists $r>0$ such that $B_r\cap \mathcal{N}=\emptyset$;
\item[(ii)] any critical point of $J_{|\mathcal{N}}$ is a critical point of $J$ on $H_0^{1,rad}(S_D\cup\Gamma_1 )$;
\item[(iii)] for any $u\in H_0^{1,rad}(S_D\cup\Gamma_1 )\setminus\{0\}$ there exists $t(u)>0$ such that $t(u)u\in\mathcal{N}$.
\end{enumerate}
\end{proposition}

Using the properties of the Nehari manifold and Proposition \ref{prop:bdd},  we also get (we refer to \cite[Theorem 3.7]{DP}):
\begin{theorem}\label{th:Nehari}
Let $f$ satisfy (F1)-(F5) with $a_1\in L^{\frac{N}{2-\epsilon}}$ for $\epsilon>0$ small enough in (F2) and assume that
$$
\frac{1}{2}f(s)s-F(s)\ge c\quad \forall s\in\R 
$$
for some constant $c\in\R$. Then there exists $\bar u \in H_0^{1}(S_D\cup\Gamma_1 )$ which is a nontrivial classical positive solution in \eqref{eq:semilinearelliptic} with Morse index equal to one.
\end{theorem}
\begin{remark} \label{rmk_subcrit}
An example of a function $f$ that satisfies the assumptions of Theorem \ref{th:Nehari} is given by 
$$
f(u)=u^p \quad 1<p<p_S.
$$
\end{remark}
\begin{remark}

We can obtain the same result by applying a mountain pass argument (\cite{DP}, Theorem 3.4).
\end{remark}

\begin{remark}\label{rmk_morse1}
It is clear that the same argument can be applied to
$$
\mathcal{N}_{rad}:=\bigg\{v\in H_0^{1,rad}(S_D\cup\Gamma_1 ):\langle J'(v),v\rangle=\int_{S_D}|\nabla v|^2dx-\int_{S_D}f(x,v)vdx=0\bigg\}.
$$
and then we obtain a solution $\tilde u\in H_0^{1,rad}(S_D\cup\Gamma_1 )$, obtained by minimization on $\mathcal{N}_{rad}$, such that $m^{rad}(\tilde u)=1$. 
\end{remark}
We conclude by proving Theorem  \ref{Th:symbreak}.
\begin{proof}[Proof of Theorem \ref{Th:symbreak}]
Let $\tilde u$ be the unique radial solution of \eqref{eq:semilinearelliptic} and let $\hat\Lambda_1^{rad}$ be the singular eigenvalue related to the linear operator $L_{\tilde  u}$.

If $\hat\Lambda_1^{rad}\ge 0$ then $m^{rad}(\tilde u)=0$. i.e. $\tilde u$ is stable in the space $H_0^{1,rad}(D\cup\Gamma_1)$. Then by $i)$ of Proposition \ref{P1} the Morse index $m(\tilde u)$ in $H_0^1(D\cup \Gamma_1)$ is also zero and hence $\tilde u$ cannot coincide with $\bar u$ which, therefore, is nonradial.

If $\hat\Lambda_1^{rad}\le 0$ and $\lambda_1(D)<-\hat\Lambda_1^{rad}$, then by $ii)$ of Proposition \ref{P1} we deduce that $m(\tilde u)\ge 2$ so that $\tilde u\ne\bar u$ and this implies the assertion, since $\tilde u$ is the only radial solution of \eqref{eq:semilinearelliptic}.
\end{proof}
\begin{remark}\label{rmksubcrit}
By Remark \ref{rmk_morse1} we know that a radial solution $\tilde u$ for which $\hat \Lambda_1^{rad}<0$ exists, under conditions $(F_1)-(F_5)$. On the other side uniqueness of the positive radial solution holds for several type of nonlinearities (see \cite{NN}). In particular if $f(u)=u^p$ $1<p<p_S$,  then there exists a unique  positive radial solution $\tilde u$ for which $m^{rad}(\tilde u)=1$.
\end{remark}
\subsection{Lane Emden }
Here we  consider the Lane-Emden problem:
\beq\label{subcritical}
\begin{cases}
-\Delta u=u^p & \textmd{ in } S_D\\
\frac{\partial u}{\partial \nu}=0 & \textmd{ on } \Gamma_1 \\
u=0& \textmd{ on } D\\
u>0 & \textmd{ in } S_D
\end{cases}
\eeq
with $1<p<p_S$,  $p_S=\frac{N+2}{N-2} N\ge 3$.

By Remarks \ref{rmk_subcrit},  \ref{rmk_morse1} and \ref{rmksubcrit}, we have that \eqref{subcritical} admits a unique positive radial solution with Morse index one in the space $H_0^{1,rad}(S_D\cup\Gamma_1)$ which is also bounded thank to Proposition \ref{prop:bdd}.   We denote it by $\tilde u_p$ to emphasize the dependence on the exponent $p$. 

If $\lambda_1(D)>N-1$, from Proposition \ref{autres} and $ii)$ of Proposition \ref{P1} we get that $m(\tilde u_p)=m^{rad}(\tilde u_p)=1$.
Indeed
$$
m(\tilde u_p)= \# 	\{j\ge 1:\lambda_j< -\hat\Lambda^{rad}_1\}+1\le \{j\ge 1:\lambda_j< N-1\}+1=1  \,.
$$
Thinking about the proof of Theorem \ref{Th:symbreak} this sems to indicate that there is no breaking of symmetry whenever $\lambda_1(D)>N-1$.

To prove instead that breaking of symmetry occurs whenever $\lambda_1(D)<N-1$, we need a  refined estimate on the eigenvalue $\Lambda_1^{rad}$ related to the linearized operator $L_{\tilde u_p}$. To stress  the dependence on $p$ we denote by $\hat \Lambda_1^{rad}(p)$ the singular radial eigenvalue.

We have
\begin{theorem}\label{eigsubcr}
Let $u\in H_0^{1,rad}(S_D\cup\Gamma_1)$ be a radial solution to \eqref{subcritical} and $\hat\Lambda_1^{rad}(p)$ be the first radial singular eigenvalue and $\psi_1$ an associated eigenfunction, i.e.
\begin{equation}\label{def:singularproblemsubc}
\begin{cases}
-\Delta \psi_1-pu^{p-1}\psi_1=\frac{\hat\Lambda^{rad}_1(p)}{|x|^2}\psi_1 & \textmd{ in } S_D\setminus\{0\}\\
\frac{\partial \psi_1}{\partial \nu}=0 & \textmd{ on }  \Gamma_1 \\
\psi_1=0& \textmd{ on }  D.
\end{cases}
\end{equation}
Then it holds 
\begin{equation} \label{lambdalimit}
\lim_{p\to p_S}\hat\Lambda^{rad}_1(p)=-(N-1).
\end{equation}
\end{theorem}
Before proving Theorem \ref{eigsubcr}, we introduce a limit eigenvalue problem (see also \cite[Section 5.1]{DIP}). Let $D^{1,2}(\R^N)$ be the space defined as
$$
D^{1,2}(\R^N)=\{u\in L^{2^*}(\R^N):|\nabla u|\in L^2(\R^N)\}
$$
and let $D^{1,2}_{rad}(\R^N)$ be its subspace given by radial function. We define
\begin{equation} \label{Lambda*1}
 \Lambda^*:=\inf_{v\in D^{1,2}_{rad} (\R^N), v\ne 0} \dfrac{\int_{\R^N}\bigg(|\nabla v|^2-p_SU^{p_S-1}v^2\bigg)dx}{\int_{\R^N}\frac{v^2}{|x|^2}dx},
\end{equation}
where
$$
U(x):=\bigg(\frac{N(N-2)}{N(N-2)+|x|^2}\bigg)^\frac{N-2}{2}.
$$
 \begin{theorem}\label{limitcase}
Let $\Lambda^*$ be given by \eqref{Lambda*1}. Then $\Lambda^*=-(N-1)$ and it is achieved by the function
\begin{equation} \label{eta_def}
\eta(x)=\frac{|x|}{\bigg(1+\frac{|x|^2}{N(N-2)}\bigg)^{N/2}} \,,
\end{equation}
which is a solution of
$$
-\Delta\eta-V\eta=\Lambda^*\frac{\eta}{|x|^2}
$$
in $\R^N\setminus\{0\}$, where $V=p_SU^{p_S-1}$.
 \end{theorem}
\begin{proof}
This theorem was proved in \cite[Theorem 5.1]{DIP2} (see also \cite[Section 5.1]{DIP},) and for this reason the proof is omitted. 
\end{proof}

\begin{proof}[Proof of Theorem \ref{eigsubcr}]
Let $p\in (1,p_S)$ be fixed and let $u_p$ the radial solution to \eqref{subcritical}. Let $\phi_p$ be the eigenfunction associated to the first eigenvalue $\hat\Lambda^{rad}_1(p)$, namely
$$
\begin{cases}
-\phi_p''-\frac{N-1}{r}\phi_p'-p|u_p|^{p-1}\phi_p=\hat\Lambda_1^{rad}(p)\frac{\phi_p}{r^2}, & r\in (0,1)\\
\phi_p(0)=\phi_p(1)=0 \,,
\end{cases}
$$
and we assume that $\phi_p$ is normalized such that
$$
\bigg\|\frac{\phi_p}{|y|}\bigg\|_{L^2(B_1)}=1.
$$
Let
$$
M_p:=\|u_p(x)\|_{L^\infty(B_{1})}<\infty
$$
and, by setting $R_p=M_p^\frac{p-1}{2}$, we rescale as follows 
$$
B_{R_p}=M_{p}^\frac{p-1}{2}B_1,
$$
and
$$
\hat \phi_p(x):=\frac{1}{M_p^\frac{(p-1)(N-2)}{4}}\phi_p\bigg(\frac{|x|}{M_p^\frac{p-1}{2}}\bigg),\quad x\in B_{R_p} \,.
$$
We have
\beq\label{eq1}
\begin{cases}
-\Delta\hat\phi_p-V_p(x)\hat\phi_p=\hat\Lambda^{rad}_1(p)\frac{\hat\phi_p}{|x|^2},  & \textmd{ in } B_{R_p}\\
\hat \phi_p= 0& \textmd{ on }\partial B_{R_p}\\
\hat \phi_p(0)= 0 \,,
\end{cases}
\eeq
where
$$
V_p(x):=p\left(\frac{1}{M_p}u_p\left(\frac{|x|}{M_p^\frac{p-1}{2}}\right) \right)^{p-1} \,,
$$
and $R_p\to\infty$ and $V_p\to V$ in $C_{loc}^0(\R^N)$ as $p\to p_S$.

Now we show \eqref{lambdalimit}. We first notice that, from Proposition \ref{autres}, we have that 
$$
\hat\Lambda^{rad}_1(p)> -(N-1)
$$
for any $1<p<p_S$, which implies 
\begin{equation}\label{B}
\liminf_{p\to p_s}\hat\Lambda^{rad}_1(p)\ge -(N-1) \,.
\end{equation}

Now we show the reverse inequality in \eqref{B}. More precisely, we show that for any $\epsilon>0$ there exists $p_\epsilon$ such that for any $p_\epsilon \leq p < p_S$ we have
$$
\hat\Lambda^{rad}_1(p)\le -(N-1)+\epsilon.
$$
Let $\epsilon>0$ be fixed. From Theorem \ref{limitcase} we know that $\Lambda^*=-(N-1)$ and it is achieved by $\eta$, which is given by \eqref{eta_def}, i.e.
$$
-(N-1)=\frac{\int_{\R^N}\bigg(|\nabla \eta|^2-p_SU^{p_S-1}\eta^2\bigg)dx}{\int_{\R^N}\frac{\eta^2}{|x|^2}dx} \,.
$$
Let
$$
\hat\eta(x)=
\begin{cases}
\eta(|x|)& |x|\in [0,R_p)\\
\eta(R_p) + \eta'(R_p)(|x|-R_p) & [R_p, \bar R_p)\\
0&[\bar R_p,\infty) \,,
\end{cases}
$$
with
$$
\bar R_p = -\eta(R_p)/\eta'(R_p) = R_p \frac{R_p^2 + N(N-2)}{(N-1)R_p^2 - N(N-2)} \,,
$$
where by an abuse of notation we write $\eta(R_p)$ for $\eta(x)$ evaluated at $x$ with $|x|=R_p$. Observe that, since $\eta$ is definitively radially monotone decreasing and convex, then there exists $p_0$ such that $\hat\eta\le \eta$ for $p\ge p_0$. 

Since $\hat \eta \in H^1_{0,rad}(B_{R_p})$, then 
$$
\hat\Lambda^{rad}_1(p)\le\frac{ \int_{B_{R_p}}(|\nabla\hat\eta(x)|^2-V_p(x)\hat\eta(x)^2)dx}{\int_{B_{R_p}} \frac{\hat\eta^2}{|x|^2}dx}
$$
and hence
$$
\hat\Lambda^{rad}_1(p)\le \frac{N_1}{D}+\frac{N_2}{D} \,,
$$
where we set
$$
N_1= \int_{B_{R_p}}(|\nabla\hat\eta(x)|^2-V(x)\hat\eta(x)^2) dx \,, \quad 
N_2= \int_{B_{R_p}}(V(x)-V_p(x))\hat\eta(x)^2 dx, \quad 
D= \int_{B_{R_p}} \frac{\hat\eta^2}{|x|^2}dx \,.
$$ 
Since $\hat\eta\le \eta$ for $p\ge p_0$ we have
\begin{equation} \label{Dbound}
D= \int_{B_{R_p}} \frac{\hat\eta^2}{|x|^2}dx\le \int_{\R^N} \frac{ \eta^2}{|x|^2}dx
\end{equation}
for $p\ge p_0$. In order to estimate $N_1$, we write
\begin{flalign*}
N_1&=\int_{B_{R_{p}}}(|\nabla\hat\eta(x)|^2-V(x)\hat\eta(x)^2)dx+\int_{B_{\bar R_p}\setminus B_{R_{p}}}(|\nabla\hat\eta(x)|^2-V(x)\hat\eta(x)^2) dx
\\&
=\int_{B_{R_{p}}}(|\nabla\eta(x)|^2-V(x)\eta(x)^2)dx+\int_{B_{\bar R_p}\setminus B_{R_{p}}}(|\nabla\hat\eta(x)|^2-V(x)\hat\eta(x)^2) dx \,.
\end{flalign*}
Since $R_p\to + \infty$ as $p \to p_S$ and $\eta \in H^1(\R^N)$, we have that for any $\epsilon>0$ there exists $p_1$ such that 
$$
\int_{B_{R_{p}}}(|\nabla\eta(x)|^2-V(x)\eta(x)^2)dx\le \int_{\R^N}(|\nabla\eta(x)|^2-V(x)\eta(x)^2)dx+\epsilon
$$
for any $p \ge p_1$. Moreover we have that 
\begin{multline*}
\Big{|} \int_{B_{\bar R_p}\setminus B_{R_{p}}}(|\nabla\hat\eta(x)|^2-V(x)\hat\eta(x)^2) dx \Big{|} \\ 
\le\int_{B_{\bar R_p}\setminus B_{R_{p}}}| \nabla\hat\eta(x)|^2 dx + \|V\|_\infty   \int_{B_{\bar R_p}\setminus B_{R_{p}}} \hat\eta(x)^2 dx 
\leq ( |\eta'(R_p)|^2  +  \|V\|_\infty  \eta(R_p)^2) |B_{\bar R_p}\setminus B_{R_{p}}| \,,
\end{multline*}
and hence there exists $p_2< p_S$ such that 
$$
\Big{|} \int_{B_{\bar R_p}\setminus B_{R_{p}}}(|\nabla\hat\eta(x)|^2-V(x)\hat\eta(x)^2) dx \Big{|}  \leq \epsilon
$$
for any $p \in [p_2,p_S)$. Thus we obtain that
$$
N_1\le  \int_{\R^N}(|\nabla\eta(x)|^2-V(x)\eta(x)^2) dx+\epsilon
$$
for any $p \geq \max(p_0,p_1,p_2)$. Since the integral on the right hand side of the above inequality is negative, from \eqref{Dbound} we have
\begin{equation}\label{N_1final}
N_1\le \frac{\int_{\R^N}(|\nabla\eta(x)|^2
-V(x)\eta(x)^2)dx }{\int_{\R^N} \frac{\eta^2}{|x|^2}dx}+\epsilon= -(N-1)+\epsilon.
\end{equation}
Now we consider $N_2$.  Let
$$
R\ge\max\bigg\{1,N(N-2), \frac{N\sqrt{(N-2)(N+2)}}{\sqrt{\epsilon}}\bigg\} 
$$
be fixed. We have
\begin{flalign*}
N_2=&\int_{B_{R_p}}\bigg[V(x)-V_{p}(x)\bigg]\hat\eta^2(x)dx\le \\
&\le\int_{B_{R_p}\cap\{|x|\le R\}}\bigg[V(x)-V_{p}(x)\bigg]\hat\eta^2(x)dx+ \int_{B_{R_p}\cap\{|x|> R\}}V(x)\hat\eta^2(x)dx\\
&=I+I\!I
\end{flalign*}
From
$$
|I|\leq \int_{B_{p}\cap\{|x|\le R\}}\big{|}V(x)-V_{p}(x)\big{|} |x|^2\frac{\hat\eta^2(x)}{|x|^2}dx\le R^2\bigg( \sup_{B_R(0)}|V_{p}(x)-V(x)|\bigg)\int_{B_{p}\cap\{|x|\le R\}} \frac{\hat \eta^2(x)}{|x|^2},
$$
since $R$ is fixed and $V_{p}\to V$ in $C^0_{loc}(\R^N)$ as $p\to p_S$, it follows that
$$
I\le \epsilon\int_{B_{p}\cap\{|x|\le R\}} \frac{\hat \eta^2(x)}{|x|^2}
$$
for $p$ close to $p_S$. 

Now we give a bound on $I\!I$. Since the function $|x|\mapsto V(x)|x|^2$ is decreasing for $|x|>N(N-2),$ our choice of $R$ yields
$$
\sup_{\{|x|>R\}}(V(x)|x|^2)\le V(R)R^2\le \frac{N^2(N+2)(N-2)}{R^2} < \epsilon \,,
$$
and then
$$
|I\!I|=\int_{B_{p}\cap\{|x|> R\}}V(x)|x|^2\frac{\hat\eta(x)^2}{|x|^2}dx\le \sup_{\{|x|>R\}}(V(x)|x|^2)\int_{B_{p}\cap\{|x|> R\}}\frac{\hat\eta^2(x)}{|x|^2}dx\le\epsilon\int_{B_{p}\cap\{|x|> R\}}\frac{\hat\eta^2(x)}{|x|^2}dx \,.
$$
Hence we have proved that 
$$
|N_2|\le  \epsilon\int_{B_{p}\cap\{|x|\le R\}} \frac{\hat \eta^2(x)}{|x|^2}\,dx+\epsilon\int_{B_{p}\cap\{|x|> R\}}\frac{\hat\eta^2(x)}{|x|^2}dx=\epsilon \,,
$$
i.e.
$$
\frac{|N_2|}{D}\le \epsilon \,,
$$
which, together \eqref{N_1final}, implies
\begin{equation} \label{questo}
\hat\Lambda^{rad}_1(p)\le -(N-1)+\epsilon \,.
\end{equation}
The assertion follows from \eqref{questo} and \eqref{B}.
\end{proof}
We are now ready to prove Theorem \ref{Th:existence}
\begin{proof}[ Proof of Theorem \ref{Th:existence}]
If $\lambda_1(D)<N-1$, by Theorem \ref{limitcase} we have that there exists $p_0\in(1,p_S)$ such that for every $p\in(p_0,p_S)$:
\beq
N-1>-\hat\Lambda_1^{rad}(p)>\lambda_1(D)
\eeq
Then by $ii)$ of Proposition \ref{P1} we have that $m(\tilde u_p)\ge 2$ for any $p\in(p_0,p_S)$.

On the other side by Theorem \ref{th:Nehari} and Remark \ref{rmk_subcrit} we know that a solution $\bar u_p$ with $m(\bar u_p)=1$ shall exist. By the uniqueness of the radial solution we get that $\bar u_p$ is not radial $\forall p\in (p_0,p_S)$. 
\end{proof}

\begin{remark}
More generally it holds that if $\lambda_j(D)<N-1$ for some $j\ge 1$ then there exists $p_0\in(1,p_S)$ such that for every $p\in(p_0,p_S),$ $
 m(\tilde u_p)\ge j+1 \,.$
\end{remark}

\section{The critical Neumann problem in the cone}
In this section we prove a symmetry breaking result for the critical Laplace equation in a  cone.  As before, we denote by $\Sigma_D\subset \R^N$, $N \geq 3$, the cone spanned by a domain $D \subset\S^{N-1}$ and consider the problem
\begin{equation}\label{def:probcone}
\begin{cases}
-\Delta u=u^{p_S}& \textmd{ in } \Sigma_D\\
\frac{\partial u}{\partial \nu}=0 & \textmd{ on }  \partial\Sigma_D\\
u>0 & \textmd{ in } \Sigma_D \,,
\end{cases}
\end{equation}
where we recall that $p_S=\frac{N+2}{N-2}=2^*-1.$  We mention that $\Sigma_D$ can be written in polar coordinates $(r,\theta)$ as  
$$
\Sigma_D = (0,+\infty) \times D \,,
$$
and in cartesian coordinates as
$$
\Sigma_D = \R^k \times \mathcal C \,,
$$
where $\mathcal C \subset \R^{N-k}$ is a cone centered at $\mathcal O \in \R^{N-k}$ which does not contains straight lines, with $k \in \{0,1,\ldots,N\}$.

Let us define the space
$$
\mathcal D^{1,2}(\Sigma_D)=\{u\in L^{2^*}:|\nabla u|\in L^2(\Sigma_D)\} 
$$
and the Sobolev quotient 
\beq\label{Sobquot}
Q_{\Sigma_D}(u)=\frac{\int_{\Sigma_D} |\nabla u|^2dx}{\bigg(\int_{\Sigma_D}|u|^{2^*}dx\bigg)^\frac{2}{2^*}} \qquad u\in D^{1,2}(\Sigma_D),  \quad u\ne 0\,.
\eeq
The following result has been proved in \cite{CP}.
\begin{theorem}\label{inf}
If $\Sigma_D$ has a point of convexity \cite[Def. 2.6]{CP} then the infimum:
\beq\label{infimum}
S_{\Sigma_D}=\inf_{u\in \mathcal{D}^{1,2}(\Sigma_D)\setminus\{0\}}Q_{\Sigma_D}(u)
\eeq
is achieved. In particular, if $\bar D\subset\S_+^{N-1}$ then $\Sigma_D$ has a point of convexity,  so that \eqref{infimum} is achieved.
\end{theorem}
From now on we consider  cones satisfying the hypothesis of Theorem \ref{inf} and observe that, since $S_{\Sigma_D}$ is achieved, the classical Sobolev inequality holds in $\mathcal{D}^{1,2}(\Sigma_D)$:
\beq\label{Sobineq}
\|u\|_{L^{2^*}(\Sigma_D)}\le(S_{\Sigma_D})^{-1/2}\|\nabla u\|_{L^2(\Sigma_D)}
\eeq
and $(S_{\Sigma_D})^{-1/2}$ is the best constant for \eqref{Sobineq}.

Then in $\mathcal{D}^{1,2}(\Sigma_D)$ we define the norm
$$
\|u\|_{\mathcal D^{1,2}(\Sigma_D)}^2=\int_{\Sigma_D}|\nabla u|^2 dx \,,
$$
which makes it a Hilbert space.

Next we consider the functional
$$
J_{\Sigma_D}=\frac{1}{2}\int_{\Sigma_D}|\nabla u|^2dx-\frac{1}{2^*}\int_{\Sigma_D} |u|^{2^*}dx,\quad u\in \mathcal D^{1,2}(\Sigma_D)
$$
together with the associated Nehari manifold 
$$
\mathcal{N} (\Sigma_D)=\bigg \{u\in \mathcal D^{1,2}(\Sigma_D):u\ne 0, \int_{\Sigma_D} |\nabla u|^2dx=\int_{\Sigma_D}|u|^{2^*}dx\bigg\}.
$$
Note that for $u\in \mathcal D^{1,2}(\Sigma_D)\setminus\{0\}$ there exists $t_u>0$ such that $t_uu\in\mathcal{N}(u)$. Then
\begin{equation}
J_{\Sigma_D}(t_u u)=\frac{1}{N}\bigg[Q_{\Sigma_D}(u)\bigg]^{\frac{N}{2}} \,,
\end{equation}
and it is easy to check that
\begin{equation} \label{inf_Sob}
C_{\Sigma_D}=\inf_{u\in \mathcal{N}_{\Sigma_D}}J_{\Sigma_D}(u)=\frac{1}{N}\bigg[S_{\Sigma_D}(u)\bigg]^\frac{N}{2} \,.
\end{equation}
We observe that any solution of \eqref{def:probcone} belongs to $\mathcal{N}_{\Sigma_D}$, which is a Hilbert manifold of codimension 1.  In particular we know that the radial function 
\begin{equation} \label{bubble}
U(x)= \alpha_N \bigg(\frac{1}{1+|x|^2}\bigg)^\frac{N-2}{2}, \quad x\in\bar\Sigma_D
\end{equation}
is a solution of (\ref{def:probcone}) with
$$
\alpha_N= (N(N-2))^\frac{N-2}{4} \,,
$$
as well as any its rescaling and admissible translation.  We will refer to $U$ as a standard bubble.  Indeed, one can see that 
$$
U_{\lambda,x_0} = \alpha_N \bigg(\frac{\lambda}{\lambda^2+|x-x_0|^2}\bigg)^\frac{N-2}{2}
$$
for $\lambda>0$ and $x_0 \in \R^k \times \{\mathcal O\}$, is a solution of (\ref{def:probcone}).

More generally, if $\Sigma_D$ is a convex cone, it has been proved in \cite{LPT} (see also \cite{CFR}) that the standard bubbles are the only solutions of (\ref{def:probcone}), up to scaling and admissible translation and we have
$$
 C_{\Sigma_D}=J_{\Sigma_D}(U)=\frac{\mathcal H^{N-1}(D)}{\mathcal H^{N-1}(\S^{N-1})}C_{\R^N}=\frac{\mathcal H^{N-1}(D)}{N \mathcal H^{N-1}(\S^{N-1})} \mathcal{S}^{N/2} \,,
$$
where $S$ is the best Sobolev constant in $\mathcal{D}^{1,2}(\R^N)$.
We also mention that, in the radial setting, the standard bubbles $U$ restricted to the cone are the only radial minimizers of $J_{\Sigma_D}$ and, actually, the only radial solutions of (\ref{def:probcone}).
   
The goal of this section is to prove Theorem \ref{Th:existencecone}, i.e. to give a condition on $D$ such that the standard bubbles are not the minimizers of $J_{\Sigma_{D}}$ on $\mathcal{N}_{\Sigma_D}$, which implies the existence of a non-radial least energy solution of (\ref{def:probcone}). This will be proven by showing that the Morse index of a bubble is strictly greater than one. Indeed, if a function $v\in \mathcal D^{1,2}(\Sigma_D)$ is a minimizer of $J_{\Sigma_D}$ on $\mathcal{N}_{\Sigma_D}$ then, since $\mathcal{N}_{\Sigma_D}$ is an Hilbert manifold of codimension 1, we have that $m(v)\le 1$ and actually $m(v)$ is exactly one for the equation \eqref{def:probcone}(see Remark \ref{rmk1}). Therefore, if we prove that $m(U)\ge 2$ for some cone $\Sigma_D$, then $U$ cannot be a minimizer of $J_{\Sigma_D}$ on $\mathcal{N}_{\Sigma_D}$ and symmetry breaking occurs.

Now recall the definition of Morse index of a solution $v$ of (\ref{def:probcone}), which is a critical point of the functional $J_{\Sigma_D}$. Since we are considering unbounded domains, some minor changes are needed with respect to the previous sections. 
We also notice that positive solutions in $\mathcal{D}^{1,2}(\Sigma_D) $ are boundd as already mentioned in Remark \ref{lastrmk}.

Let $Q_v(\cdot)$ be the quadratic form corresponding to the linearization of the equation \eqref{def:probcone} at a solution $v$
$$
 Q_v(\psi)=\int_{\Sigma_D}|\nabla\psi|^2 dx-p_S\int_{\Sigma_D} |v|^{p_S-1}\psi^2 dx\qquad \forall\psi\in C_c^1(\bar\Sigma_D)
$$
and notice that $Q_v$ it the quadratic form corresponding to the second derivative of the functional $J_{\Sigma_D}$ in $\mathcal D^{1,2}(\Sigma_D)$. 
 \begin{definition}
 Let $v\in \mathcal D^{1,2}(\Sigma_D)$be a solution of (\ref{def:probcone}). We say that:
 \begin{enumerate}
 \item[(i)] $v$ is stable (or has zero Morse index) if $Q_v(\psi)\ge 0$ for any $\psi\in C_c^1(\bar\Sigma_D)$;
 \item[(ii)] $v$ has Morse index equal to the integer $m(v)\ge 1$ if $m(v)$ is the maximal dimension of a subspace of $C^1_c(\bar \Sigma_D)$ where the quadratic form $Q_v$ is negative semidefinite;
 \item[(iii)] $v$ has infinite Morse index if, for any integer $k\ge 1$, there exists a k-dimensional subspace of $C^1_c(\bar\Sigma_D)$ where $Q_v$ is negative definite.
 \end{enumerate}
 \end{definition}

As in Section 3 we denote by $\lambda_j(D)$, $j \in \mathbb N$, the $j-$nontrivial eigenvalue of the Laplace-Beltrami operator $-\Delta_{\S^{N-1}}$ on the domain $D$ with zero Neumann boundary condition.  We prove the following result:
  \begin{theorem}\label{Th:formula}
Let $U$ be given by \eqref{bubble}, which is a solution to \eqref{def:probcone}. We have  
\begin{equation} \label{mUthm}
m(U)=\# \{j\ge 1|\lambda_j(D)< N-1\}+1\, 
\end{equation}
where $m(U)$ is the Morse index of $U$. Therefore it holds
\begin{enumerate}
 \item[(i)] if $\lambda_1(D)<N-1$ then $m(U)\ge 2$;
 \item[(ii)] if $\lambda_1(D)\ge N-1$ then $m(U)=1$.
\end{enumerate}
\end{theorem} 

Though quite obvious it is important to remark that the same result holds if  instead of $U$ we consider any rescaling of it. This can be verified by analyzing how the quadratic form changes for the rescaled bubbles. 
Hence, if a cone satisfies (i) in Theorem \ref{Th:formula} then the standard bubble $U$ cannot achieve the infimum $C_{\Sigma_D}$ and then it cannot be a least-energy solution of (\ref{def:probcone}). 
 
 We delay the proof of Theorem \ref{Th:formula} and show how to deduce Theorem \ref{Th:existencecone} from Theorem \ref{Th:formula}.
 \begin{proof}[Proof of Theorem \ref{Th:existencecone}]
From Theorem \ref{Th:formula} we first deduce that the bubbles cannot minimize $C_{\Sigma_D}$  whenever $\lambda_1(D)<N-1.$ Moreover, since $\bar D\subset\S^{N-1}_+$ then Theorem \ref{inf} yields that $C_{\Sigma_D}$ is achieved. Since the only radial solutions of \eqref{def:probcone} are the bubbles, we get the existence of a nonradial solution $w$ of \eqref{def:probcone} which is actually a least energy solution. Moreover, by using Kelvin transform and regularity results (see \cite{CFR}), we also deduce that $w$ is fast decaying.
 \end{proof}
 
 \begin{remark}\label{rmk1}
 To prove that the quadratic form $Q_U$ becomes negative for two independent directions let us first observe that it is negative on the space spanned by $U$ itself.  Indeed, since $p_S=2^*-1$, from the equation (\ref{def:probcone}) we get:
$$
 \int_{\Sigma_D} |\nabla U|^2dx=\int_{\Sigma_D} U^{2^*}dx
$$
and then
$$
Q_U(U)=\int_{\Sigma_D}(U^{2^*}-p_SU^{2^*})dx<0 
$$
since $p_S>1$.
 \end{remark}

 Let $\mathcal D_{rad}^{1,2}(\R^N)$ and $\mathcal D^{1,2}_{rad}(\Sigma_D)$ be, respectively, the subspaces of $\mathcal D^{1,2}(\R^N)$ and $\mathcal D^{1,2}(\Sigma_D)$ given by radial functions. We observe that, since $N\ge 3$, Hardy inequality in Proposition \ref{prop_Hardy} yields that if $v \in \mathcal D^{1,2}(\R^N)$ and $w \in \mathcal D^{1,2}(\Sigma_D)$ then $\frac{v}{|x|}\in L^2(\R^N)$ and $\frac{w}{|x|}\in L^2(\Sigma_D)$. Hence, as observed in Section 3, the eigenvalue $\Lambda^*$ given by \eqref{Lambda*1} is well defined and Theorem \ref{limitcase} holds. 
 We conclude by proving Theorem \ref{Th:formula}.

 \begin{proof}[Proof of Theorem \ref{Th:formula}]
We start by proving \eqref{mUthm}. This will be done in two steps.
 
\textbf{Step 1.} We prove that
$$
m(U)\ge\# \{j\ge 1|\lambda_j< N-1\}+1.
$$
 Suppose that $\# \{j\ge 1|\lambda_j< N-1\}=k$, we want to prove that $m(U)\ge k+1.$ From Remark \ref{rmk1}, we know that $Q_U(U)< 0$, which implies $m(U)\ge 1$. Now we want to find $\psi_1,...,\psi_k$ non radial and linearly independent, such that $Q_U(\psi_i)<0,$ for $i=1,...,k.$
Let $i\in\{1,...,k\}$ be fixed and assume that $\lambda_i(D)<N-1$.
Let $\varphi_i$ be an eigenfunction of $-\Delta_{\S^{n-1}}$ associated to $\lambda_i(D)$. Hence $\varphi_i$ is a solution of the Neumann problem
$$
 \begin{cases}
 -\Delta_{\S^{n-1}}\varphi_i=\lambda_i(D)\varphi_i&\text{in }D\\
 \frac{\partial\varphi_i}{\partial\nu}=0 &\text{on }\partial D.
 \end{cases}
$$
Now we consider the function
$$
 \psi_i(x)=\eta(|x|)\varphi_i\bigg(\frac{x}{|x|}\bigg)\quad x\in \Sigma_D \,,
$$
where $\eta$ is defined in \eqref{eta_def} and notice that we have (for $r=|x|$)
 \begin{flalign*}
 -\Delta\psi_i&=(-\eta''-\frac{n-1}{r}\eta'_1)\varphi_i-\frac{\eta}{r^2}\Delta_{\S^{n-1}}\varphi_i=\\
 &=\left [p_S|U|^{p_S-1}\eta+\frac{\Lambda^*}{r^2}\eta\right ]\varphi_i+\frac{\lambda_i(D)}{r^2}\eta\varphi_i=\\
 &=p_S|U|^{p_S-1}\psi_i+\frac{\Lambda^*+\lambda_i(D)}{r^2}\psi_i 
 \end{flalign*}
in $\Sigma_D$. Hence
$$
 -\Delta\psi_i-p_S|U|^{p_S-1}\psi_i=(\Lambda^*+\lambda_1(D))\frac{\psi_i}{r^2}.
$$
By multiplying by $\psi_i$, integrating by parts and using Neumann condition for $\eta$ and $\varphi_i$, we get that
$$
 Q_U(\psi_i)<0 \,.
$$
It is clear that $\psi_i$ is not radial and $\psi_1,...,\psi_k$ are linearly independent because $\varphi_1,..,\varphi_k$ are linearly independent.

\medskip

  \textbf{Step 2} We prove that
  $$
m(U)\le\# \{j\ge 1|\lambda_j< N-1\}+1.
$$
If $m(U)=k$ there exist $\varphi_1,...,\varphi_{k}\in C_c^1(\bar\Sigma_D)$, linearly independent, s.t.  
  $$
  Q_U(\varphi_i,\varphi_i)<0
  $$
for any $i=1,...,k$. We consider a ball $B_R$ centered at the origin and of radius $R$ such that
  $$
  {\rm supp} \{\varphi_1,...,\varphi_{k}\}\subset B_R.
  $$
Hence, for any $i=1,...,k$, we have
  $$
  Q_U(\varphi_i)=\int_{\Sigma_D}|\nabla\varphi_i|^2 dx-p_S\int_{\Sigma_D} |U|^{p_S-1}\varphi_i^2 dx=\int_{\Sigma_D\cap B_R}|\nabla\varphi_i|^2 dx-p_S\int_{\Sigma_D\cap B_R} |U|^{p_S-1}\varphi_i^2 dx \,.
  $$
Let $\bar u_{\epsilon,R}\in H_0^{1,rad}(\Sigma_D)$ be a radial solution of the subcritical problem
  $$
\begin{cases}
-\Delta u=u^{p_S-\epsilon}&  \textmd{ in } \Sigma_D\cap B_R\\
u=0& \textmd{ on } \Sigma_D\cap \partial B_R\\
\partial_\nu u=0& \textmd{ on } \partial\Sigma\cap \bar B_R\\
u>0& \textmd{ in } \Sigma_D\cap B_R \,,
\end{cases}
  $$
and consider the associated quadratic form
   $$
  Q_{\bar u_{\epsilon,R}}(\varphi)=\int_{ \Sigma_D\cap B_R}|\nabla\varphi|^2 dx-(p_S-\epsilon)\int_{ \Sigma_D\cap B_R} |\bar u_{\epsilon,R}|^{p_S-\epsilon-1}\varphi^2 dx.
  $$
Let $M_\epsilon$ be the maximum of $\bar u_{\epsilon,R}$ (which is known to be attained at the origin)
 $$
 M_\epsilon:=\bar u_{\epsilon,R}(0)=\|\bar u_{\epsilon,R}\|_{L^\infty}\to \infty \quad\text{ for }\epsilon\to 0
 $$
 and we define
 $$
 \tilde{u}_\epsilon(x):=\frac{1}{M_\epsilon}\bar u_{\epsilon,R}\bigg(\frac{x}{M_\epsilon^{\frac{p_\epsilon-1}{2}}}\bigg).
 $$
 If we call $B_{R_\epsilon}=M^{\frac{p-1}{2}}_\epsilon B(0,R)$ we have that $\tilde{u}_\epsilon$ satisfies
 $$
 \begin{cases}
 -\Delta\tilde{u}_\epsilon=\tilde{u}_\epsilon^{p_\epsilon} &  \textmd{ in }B_{R_\epsilon}\\
 \tilde{u}_\epsilon=0&   \textmd{ on } \partial B_{R_\epsilon}\\
\partial_\nu\tilde u_\epsilon=0& \textmd{ on }  \partial B_{R_\epsilon}\\
\tilde u_\epsilon>0 & \textmd{ in } B_{R_\epsilon}\\
\tilde{u}_\epsilon(0)=1 \,.
 \end{cases}
 $$
By a standard compactness argument \cite[Section 5.1]{DIP} we have that  
\begin{equation} \label{utoU}
 \tilde u_\epsilon\to U \quad\text{for } \epsilon\to 0 \quad\text{ in }C^2_{loc}(\R^N)
\end{equation}
where $U$ is the standard bubble with $U(0)=1$.
 
If $\epsilon$ is small enough we have that $M_\epsilon^\frac{p_\epsilon-1}{2}>R$. Hence, for any fixed $i=1,...,k$, from \eqref{utoU} we have that 
$$
  Q_{\tilde u_{\epsilon}}(\varphi_i)=\int_{\text{supp}(\varphi_i)}|\nabla\varphi_i|^2 dx-(p_S-\epsilon)\int_{\text{supp}(\varphi_i)} |\tilde u_{\epsilon}|^{p_S-\epsilon-1}\varphi_i^2 dx\to Q_U(\varphi_i)<0\,,
$$
as $\epsilon\to 0$.
It follows that $m(\tilde u_{\epsilon}) \geq k$; since $ii)$ of Proposition \ref{P1} and Theorem \ref{eigsubcr} imply
$$
m(\tilde u_{\epsilon}) = \# \{j\ge 1|\lambda_j< N-1\}+1
$$
we conclude the proof of Step 2.
\end{proof}

\end{document}